\overfullrule=0pt
\input amstex
\documentstyle{amsppt}
\newcount\mgnf\newcount\tipi\newcount\tipoformule\newcount\greco
\tipi=2          
\tipoformule=0   

\global\newcount\numsec\global\newcount\numfor
\global\newcount\numapp\global\newcount\numcap
\global\newcount\numfig\global\newcount\numpag
\global\newcount\numnf
\global\newcount\numtheo

\def\SIA #1,#2,#3 {\senondefinito{#1#2}%
\expandafter\xdef\csname #1#2\endcsname{#3}\else
\write16{???? ma #1,#2 e' gia' stato definito !!!!} \fi}

\def \FU(#1)#2{\SIA fu,#1,#2 }

\def\etichetta(#1){(\veroparagrafo.\veraformula)%
\SIA e,#1,(\veroparagrafo.\veraformula) %
\global\advance\numfor by 1%
\write15{\string\FU (#1){\equ(#1)}}%
\write16{ EQ #1 ==> \equ(#1)  }}

\def\etichettat(#1){\veroparagrafo.\veratheorema:%
\SIA e,#1,{\veroparagrafo.\veratheorema} %
\global\advance\numtheo by 1%
\write15{\string\FU (#1){\thu(#1)}}%
\write16{ TH #1 ==> \thu(#1)  }}

\def\etichettaa(#1){(A\veraappendice.\veraformula)
 \SIA e,#1,(A\veraappendice.\veraformula)
 \global\advance\numfor by 1
 \write15{\string\FU (#1){\equ(#1)}}
 \write16{ EQ #1 ==> \equ(#1) }}
\def\getichetta(#1){Fig. \verafigura
 \SIA g,#1,{\verafigura}
 \global\advance\numfig by 1
 \write15{\string\FU (#1){\graf(#1)}}
 \write16{ Fig. #1 ==> \graf(#1) }}
\def\retichetta(#1){\numpag=\pgn\SIA r,#1,{\verapagina}
 \write15{\string\FU (#1){\rif(#1)}}
 \write16{\rif(#1) ha simbolo  #1  }}
\def\etichettan(#1){(n\verocapitolo.\veranformula)
 \SIA e,#1,(n\verocapitolo.\veranformula)
 \global\advance\numnf by 1
\write16{\equ(#1) <= #1  }}

\newdimen\gwidth
\gdef\profonditastruttura{\dp\strutbox}
\def\senondefinito#1{\expandafter\ifx\csname#1\endcsname\relax}
\def\BOZZA{
\def\alato(##1){
 {\vtop to \profonditastruttura{\baselineskip
 \profonditastruttura\vss
 \rlap{\kern-\hsize\kern-1.2truecm{$\scriptstyle##1$}}}}}
\def\galato(##1){ \gwidth=\hsize \divide\gwidth by 2
 {\vtop to \profonditastruttura{\baselineskip
 \profonditastruttura\vss
 \rlap{\kern-\gwidth\kern-1.2truecm{$\scriptstyle##1$}}}}}
\def\verapagina{
{\romannumeral\number\numcap}.\number\numsec.\number\numpag}}

\def\alato(#1){}
\def\galato(#1){}
\def\veroparagrafo{\number\numsec}\def\veraformula{\number\numfor}
\def\veraappendice{\number\numapp}
\def\verapagina{\number\pageno}\def\veranformula{\number\numnf}
\def\verafigura{{\romannumeral\number\numcap}.\number\numfig}
\def\verocapitolo{\number\numcap}\def\veranformula{\number\numnf}
\def\veratheorema{\number\numtheo}
\def\Eqn(#1){\eqno{\etichettan(#1)\alato(#1)}}
\def\eqn(#1){\etichettan(#1)\alato(#1)}
\def\TH(#1){{\etichettat(#1)\alato(#1)}}
\def\thv(#1){\senondefinito{fu#1}$\clubsuit$#1\else\csname fu#1\endcsname\fi}
\def\thu(#1){\senondefinito{e#1}\thv(#1)\else\csname e#1\endcsname\fi}

\def\Eq(#1){\eqno{\etichetta(#1)\alato(#1)}}
\def\eq(#1){\etichetta(#1)\alato(#1)}
\def\Eqa(#1){\eqno{\etichettaa(#1)\alato(#1)}}
\def\eqa(#1){\etichettaa(#1)\alato(#1)}
\def\dgraf(#1){\getichetta(#1)\galato(#1)}
\def\drif(#1){\retichetta(#1)}

\def\eqv(#1){\senondefinito{fu#1}$\clubsuit$#1\else\csname fu#1\endcsname\fi}
\def\equ(#1){\senondefinito{e#1}\eqv(#1)\else\csname e#1\endcsname\fi}
\def\graf(#1){\senondefinito{g#1}\eqv(#1)\else\csname g#1\endcsname\fi}
\def\rif(#1){\senondefinito{r#1}\eqv(#1)\else\csname r#1\endcsname\fi}
\def\bib[#1]{[#1]\numpag=\pgn
\write13{\string[#1],\verapagina}}

\def\include#1{
\openin13=#1.aux \ifeof13 \relax \else
\input #1.aux \closein13 \fi}

\openin14=\jobname.aux \ifeof14 \relax \else
\input \jobname.aux \closein14 \fi
\openout15=\jobname.aux
\openout13=\jobname.bib


\ifnum\tipoformule=1\let\Eq=\eqno\def\eq{}\let\Eqa=\eqno\def\eqa{}
\def\equ{}\fi


{\count255=\time\divide\count255 by 60 \xdef\hourmin{\number\count255}
        \multiply\count255 by-60\advance\count255 by\time
   \xdef\hourmin{\hourmin:\ifnum\count255<10 0\fi\the\count255}}

\def\oramin{\hourmin }

\def\data{\number\day/\ifcase\month\or january \or february \or march \or
april \or may \or june \or july \or august \or september
\or october \or november \or december \fi/\number\year;\ \oramin}

\def\titdate{ \ifcase\month\or January \or February \or March \or
April \or May \or June \or July \or August \or September
\or October \or November \or December \fi \number\day, \number\year;\ \oramin}

\def\titdatebis{ \ifcase\month\or January \or February \or March \or
April \or May \or June \or July \or August \or September
\or October \or November \or December \fi \number\day, \number\year}


\newcount\pgn \pgn=1
\def\foglio{\number\numsec:\number\pgn
\global\advance\pgn by 1}
\def\foglioa{A\number\numsec:\number\pgn
\global\advance\pgn by 1}

\footline={\rlap{\hbox{\copy200}}\hss\tenrm\folio\hss}


\global\newcount\numpunt

\magnification=\magstephalf
\baselineskip=12pt
\parskip=6pt

\hoffset=1.0truepc 
\hsize=6.1truein
\vsize=8.4truein 

\def\a{\alpha}
\predefine\barunder{\b}
\redefine\b{\beta}
\def\d{\delta}
\def\e{\epsilon}

\def\f{\phi}
\def\g{\gamma}

\def\l{\lambda}

\def\s{\sigma}
\def\t{\tau}

\def\o{\omega}
\def\D{\Delta}
\def\L{\Lambda}
\def\G{\Gamma}
\def\O{\Omega}

\def\[{{[\kern-.18em{[}}}
\def\]{{]\kern-.18em{]}}}

\def\1{{1\kern-.25em\roman{I}}}
\def\eu{{1\kern-.25em\roman{I}}}
\def\f1{{1\kern-.25em\roman{I}}}

\def\R{{\Bbb R}}  
\def\N{{\Bbb N}}  
\def\P{{\Bbb P}}  
\def\Q{{\Bbb Q}}  
\def\E{{\Bbb E}}  

\def\dist{\,\roman{dist}}

\let\cal=\Cal

\def\CC{{\cal C}}

\def\EE{{\cal E}}
\def\FF{{\cal F}}
\def\GG{{\cal G}}

\def\OO{{\cal O}}
\def\PP{{\cal P}}

\def\RR{{\cal R}}

\def\TT{{\cal T}}
\def\VV{{\cal V}}

\def\VV{{\cal V}}

\def\chap #1#2{\line{\ch #1\hfill}\numsec=#2\numfor=1\numtheo=1}

\def\wt{\widetilde}
\def\wh{\widehat}


\def\note#1{\footnote{#1}}

\def\frac#1#2{{#1\over #2}}
\def\sfrac#1#2{{\textstyle{#1\over #2}}}

\def\text#1{\quad{\hbox{#1}}\quad}
\def\newpage{\vfill\eject}
\def\proposition #1{\noindent{\thbf Proposition #1}}

\def\theo #1{\noindent{\thbf Theorem {#1} }}

\def\lemma #1{\noindent{\thbf Lemma {#1} }}
\def\definition #1{\noindent{\thbf Definition {#1} }}

\def\proof{{\noindent\pr Proof: }}
\def\proofof #1{{\noindent\pr Proof of #1: }}
\def\endproof{\hfill$\square$}
\def\remark{\noindent{\bf Remark: }}
\def\thanks{\noindent{\bf Acknowledgements: }}

\font\pr=cmbxsl10

\font\thbf=cmbxsl10 scaled\magstephalf

\font\ch=cmbx12

\font\it=cmti10
\font\bf=cmbx10

\def\PRM{\hbox{\rm PRM}}
\def\asl{\hbox{\rm Asl}}



\font\tit=cmbx12
\font\aut=cmbx12

\def\s{\char'31}

\centerline{\tit AGING IN REVERSIBLE DYNAMICS OF DISORDERED SYSTEMS.}
\smallskip
\centerline{\tit II. emergence of the arcsine law in the random hopping}
\centerline{\tit time dynamics of the REM.}

\vskip.2truecm
\vskip1truecm
\vskip.5cm

\centerline{\aut V\'eronique  Gayrard
\note{
 CMI, LATP, Universit\'e de Provence
 39 rue F. Joliot-Curie, 13453 Marseille cedex 13\hfill\break
e-mail: gayrard\@cmi.univ-mrs.fr}
}

\vskip0.5cm
\centerline{\titdatebis}

\vskip.5cm

\vskip0.5truecm\rm
\def\s{\sigma}
\noindent {\bf Abstract:} Applying the new tools developed in \cite{G1}, we investigate
the arcsine aging regime of
the random hopping time dynamics of the REM. Our results are optimal in several ways.
They cover the full time-scale and temperature domain where this phenomenon occurs.
On this domain the limiting clock process and associated time correlation function are
explicitly constructed. Finally, all convergence statements w.r.t\. the law of the random environment
are obtained in the strongest sense possible,
except perhaps on the very last scales before equilibrium.


\vfill
$ {} $

\newpage




\chap{1. Introduction.}1

\smallskip
\line{\bf 1.1. The setting.\hfill}

New tools for the study of the aging behavior of disordered systems were developed
in \cite{G1} and successfully applied to the study of the arcsine aging regime of
Bouchaud's asymmetric trap model on the complete graph \cite{B, BD, BRM}.
We refer the reader to the introduction of \cite{G1}
for an overall presentation of the aging phenomenon.
In this follow up paper we continue the investigation of the aging behavior
of disordered systems, focusing on the the so-called random hoping time dynamics
of Derrida's REM \cite{D1,D2}. Our objectives here are twofold:

\noindent 1) Establish the optimal time scale and temperature domain where the model
exhibits
an arcsine aging regime,
striving for
statements that are valid in the strongest sense possible w.r.t\. the law of the random environment.
Indeed, if the random hoping time dynamics (RHT) of the REM has been intensively studied,
each existing contribution \cite{BBG1, BBG2, BBG3, BC, CG} concerns
specific choices of the time scale and temperature parameters, and
a full, unifying picture was still missing.

\noindent 2) Explain why and how the RHT dynamics of the REM and Bouchaud's
REM-like trap model
exhibit the same arcsine aging regime. More precisely, identify and isolate
the mechanisms that, within the new technical framework of \cite{G1}, allow to reduce the study
of (the clock process of) the REM to that of its trap version.
This technical issue is a crucial step towards the understanding of the
more interesting metropolis dynamics of the REM, \cite{G2}.

We now specify the model and succinctly recall the basics of arcsine aging
-- a detailed exposition can be found in \cite{G1}.
Denote by $\VV_n=\{-1,1\}^n$ the n-dimensional discrete cube, and
by $\EE_n$ its edges set.
On $\VV_n$ we construct a random landscape of traps (the random environment) by assigning
to each site, $x$, the Boltzman weight of the REM, $\t_n(x)$. Namely,
given a parameter $\beta>0$, called the inverse temperature, and a collection
$(H_n(x), x\in\VV_n)$ of independent standard gaussian random variables,
we set
$
\t_n(x)=\exp(-\b \sqrt n H_n(x))
$.
(The dependence of $\t_n(x)$ on  $\beta$ will be kept implicit.)
The sequence $(\t_n(x), x\in\VV_n)$, $n\geq 1$, is
defined on a common probability space denoted $(\O^{\t}, \FF^{\t}, \P)$.

The RHT dynamics in the landscape $(\t_n(x), x\in\VV_n)$
is a continuous time Markov chain  $(X_{n}(t), t>0)$ on $\VV_n$
that
can be constructed  as follows:
let $(J_n(k), k\in\N)$ be the
simple random walk on $\VV_n$ with initial distribution $\mu_n$ and transition probabilities
$$
p_n(x,y)=\frac 1n\,,\quad\forall(x,y)\in \EE_n\,,
\Eq(1.9)
$$
and let the {\it clock process} be the partial sum process
$$
\wt S_n(k)=\sum_{i=0}^{k}\t_n(J_n(i))e_{n,i}\,,\quad k\in \N\,,
\Eq(1.1.6)
$$
where $(e_{n,i}\,,n\in\N, i\in\N)$ is a family of independent
mean one exponential random variables, independent of $J_n$;
then
$$
X_n(t)=J_n(i) \text{if} \wt S_n(i)\leq t<\wt S_n(i+1)  \text{for some} i\,.
\Eq(1.1.7)
$$
This defines $X_n$ in terms of its jump chain, $J_n$, and holding times, $\t_n(x)$ being the mean value of the
exponential holding time at $x$. Equivalently, $X_n$ is the chain with initial distribution $\mu_n$ and
jump rates
$
\l_n(x,y)=(n\t_n(x))^{-1}\,
$,
$(x,y)\in \EE_n$, $x\neq y$.
This last description makes it easy to check that $X_n$ is a Glauber dynamics, namely,
that it is reversible with respect to the measure (the Gibbs measure of the REM)
defined on $\VV_n$ by
$$
\GG_n(x)=\frac{\t_n(x)}{{\textstyle\sum_{x\in\VV_n}\t_n(x)}}\,,\quad x\in \VV_n\,.
\Eq(1.7)
$$
The model we referred to as the REM-like trap model
was proposed by Bouchaud
as an approximation of
the aging dynamics of the REM (see \cite{BBG1} for details on this derivation).
It is a Markov chain
$X'_n$ on $\VV'_n=\{1,\dots,n\}$ with jump rates
$
\l'_n(x,y)=(n\t'(x))^{-1}\,
$,
$(x,y)\in \VV'_n\times\VV'_n$, $x\neq y$, where
$(\t'(x), x\in\VV'_n)$ are i.i.d\. r.v\. in the domain
of attraction of a positive stable law with index $0<\a<1$\note{
The model we just described
is sometimes called Bouchaud's trap model on the complete graph.
It is obtained by setting $a=0$ in Bouchaud's
asymmetric trap model on the complete graph studied in \cite{G1}.
In Bouchaud's original model, the following choices are made:
$
\t'(x)=e^{E_x/\a}
$,
where $(E_x, x\in\VV'_n)$ are i.i.d\. mean one exponential r.v\.'s,
and $\a=\sqrt{2\log2}/\beta$, $\beta>0$.
}.

To study the aging behavior of $X_n$ we need to choose three ingredients:
1) an initial distribution, $\mu_n$;
2) a time scale of observation, $c_n$;
and 3) a time-time correlation function, $\CC_{n}(t,s)$, $t, s\geq 0$:
this is a function that quantifies
the correlation between the state of the system at time $t$,
$X_n(c_nt)$, and its state at time $t+s$,  $X_n(c_n(t+s))$.
%
%
We will say that the process $X_n$ has an arcsine aging regime
of parameter $0<\a<1$ whenever one can find
a time-time correlation function such that, denoting by $\asl_\a(\cdot)$
the generalized  arcsine distribution function of parameter $\a$,
$$
\asl_{\a}(u)=\frac{\sin{\a}\pi}{\pi}\int_0^u (1-x)^{-\a}x^{-1+\a}dx\,,\quad 0<\a< 1\,,
\Eq(1.theo1.0'')
$$
one of the following relations\note{
These are the only two limiting procedure relevant for the REM.
A third one was present in Bouchaud's REM-like trap model:
$
\lim_{t\rightarrow\infty}\lim_{n\rightarrow\infty}
\CC_{n}(t,\rho t)=\asl_{\a}(1/1+\rho)
$.
}
holds true,
$$
\lim_{t\rightarrow 0}\lim_{n\rightarrow\infty}
\CC_{n}(t,\rho t)=\asl_{\a}(1/1+\rho)\,,
\Eq(1.1.14)
$$
$$
\lim_{n\rightarrow\infty}
\CC_{n}(t,\rho t)=\asl_{\a}(1/1+\rho)\,,\,\,\,t>0\,\,\,\hbox{\rm arbitrary,}
\Eq(1.1.13)
$$
for all $\rho\geq 0$, and some convergence mode w.r.t\.
the probability law $\P$ of the random landscape.
 It is today well understood that the existence of
an arcsine aging regime is governed by Dynkin and Lamperti's arcsine law for subordinators, applied
to the limiting, appropriately re-scaled, clock process:
arcsine aging will be present when the re-scaled clock process converges to a subordinator whose L\'evy measure satisfies
the slow variation conditions of the Dynkin-Lamperti Theorem\note{See e.g\. Appendix A.2.1 of \cite{G1} for a statement of the latter.}.
In this light a natural
indicator of arcsine aging, and the one we will choose,
is:
$$
\CC_{n}(t,s)=\PP_{\mu_n}\left(\left\{c_n^{-1}\wt S_n(k)\,, k\in\N\right\}
\cap (t, t+s)=\emptyset\right)\,,\quad 0\leq t<t+s\,.
\Eq(1.1.8)
$$
where $\PP_{\mu_n}$ denotes the law of $X_n$ with initial distribution $\mu_n$.
In words, this is the probability that the range of the re-scaled clock process $c_n^{-1}\wt S_n$
does not intersect the time interval $(t, t+s)$.
The initial distribution is taken to be the invariant measure of the jump chain;
that is we set  $\mu_n=\pi_n$, where
$$
\pi_n(x)=2^{-n}\,,\quad x\in \VV_n\,.
\Eq(1.8)
$$
We will see in later that this choice is generic
\note{
See the remark below  Proposition \thv(4.prop1).
}.

We now proceed to state our results. Subsection 1.2 contains the results on the convergence of $\CC_{n}$.
Their parent results on the convergence of the re-scaled clock process will be stated next, in Subsection 1.3.
The closing Subsection 1.4 is devoted to the method of proof:
we recall the needed results from \cite{G1}
and elaborate on the strategy of the proofs.

\bigskip
\line{\bf 1.2. Convergence of the time-time correlation  $\CC_{n}$.\hfill}

To keep the notation of this paper consistent with those of \cite{G1} we set $c_n=r_n$ and call $r_n$
a space scale.\note{We saw in \cite{G1} (see the paragraph above
Definition 4.1) that, generally speaking, measuring time on scale $c_n$ for
$X_n$ corresponds to re-scaling the traps $\t_n(x)$ by a certain amount $r_n$.
It is clear from \eqv(1.1.7) that here space and time scales coincide.
}
We will distinguish three  types of space scales:
the {\it short} scales, the {\it intermediate} scales and the {\it extreme} scales.
Let $b_n$ be defined through
$$
b_n\P(\t_n(x)\geq r_n)=1\,,
\Eq(1.1')
$$
and set $m_n=\log b_n/\log2$.

\definition{\TH(1.def1)}
We say that a diverging sequence $r_n$ is:
\item{(i)} a {\it short space scale} if $m_n$ is a diverging sequence such that
$$
\frac{m_n}{n}=o(1)\,,
\Eq(1.1'')
$$
\item{(ii)} an {\it intermediate space scale} if there exists $0<\varepsilon\leq 1$ such that
$$
\frac{m_n}{n}\sim\varepsilon  \text{and} \frac{2^{m_n}}{2^n}=o(1)\,,
\Eq(1.2)
$$
\item{(iii)} an {\it extreme space scale} if there exists $0<\varepsilon'\leq 1$ such that
$$
\frac{2^{m_n}}{2^n}\sim \varepsilon'\,.
\Eq(1.3)
$$

Note that Definition \thv(1.def1) allows to classify all sequences $b_n$ such that $1\ll b_n\leq 2^n$,
leaving out constant sequences only.
For any of the above space scale set $\varepsilon=\lim_{n\rightarrow\infty}\frac{m_n}{n}$. Thus $\varepsilon=0$ if $r_n$ is a short space scale,
$0<\varepsilon\leq 1$ if $r_n$ is an intermediate space scale, and $\varepsilon=1$ if $r_n$ is an extreme space scale.
For  $0\leq\varepsilon\leq 1$ and  $0<\beta<\infty$, define
$$
\eqalign{
\b_c(\varepsilon)&=\sqrt{\varepsilon2\log 2}\,,
\cr
\a(\varepsilon)&=\b_c(\varepsilon)/\b\,,
}
\Eq(1.theo1.0')
$$
and write $\b_c\equiv\b_c(1)$ and $\a\equiv\a(1)$.

Before stating our two theorems let us outline their content.
Theorem \thv(1.theo1) establishes that arcsine aging is not present
on short scales, where $\varepsilon=0$, and emerges on intermediate
scales, as $\varepsilon$ becomes positive. There, it is present in
the entire time scale and temperature domain
$
\left\{0<\varepsilon\leq 1 ; 0<\a(\varepsilon)<1\right\}
$,
and gets interrupted on the critical temperature line
$
\left\{0<\varepsilon\leq 1 ; \a(\varepsilon)=1\right\}
$.
These results hold true for all time $t>0$, that is, \eqv(1.1.13) prevails.
On extreme scales,  where $\varepsilon=1$, arcsine aging is still present on
the entire low temperature line
$
\left\{\varepsilon= 1 ; 0<\a\equiv\a(1)<1\right\}
$,
but in the limit $t\rightarrow 0$ only, that is, along the double
limiting procedure of \eqv(1.1.14). These are the very last (and longest)
scales where arcsine aging occurs before interruption. Indeed, as $t$
increases from $0$ to $\infty$, the system moves out of its
arcsine aging regime and crosses over to its stationarity regime.
This last statement is the content of Theorem \thv(1.theo2).
Clearly all these convergence results must hold in some sense w.r.t\.
the law $\P$ of the random landscape. We now make this precise.

\theo{\TH(1.theo1)}
{\it 
Take $\mu_n=\pi_n$.

\item{(i)}
Let $r_n$ be a short space scale and, given any $0<\beta<\infty$, assume that
$
\frac 1{\b  n}\log r_n\geq 4\sqrt{\sfrac{\log n}{n}}
$.
Then, $\P$-almost surely, for all $t\geq 0$ and all $\rho>0$,
$$
\lim_{n\rightarrow\infty}\CC_{n}(t,\rho t)=1\,.
\Eq(1.theo1.0)
$$

\item{(ii)}
Let $r_n$ be an intermediate space scale.
For all  $0<\varepsilon\leq 1$ and all  $0<\beta<\infty$ such that
 $0<\a(\varepsilon)\leq 1$ the following holds,
$\P$-almost surely if
$\frac{2^{m_n}}{2^n}\log n=o(1)$,
and in $\P$-probability otherwise:
for all $t\geq 0$ and all $\rho>0$,
$$
\lim_{n\rightarrow\infty}\CC_{n}(t,\rho t)=
\cases
\asl_{\a(\varepsilon)}(1/1+\rho) , &\hbox{if $\a(\varepsilon)< 1$,}\,\,\,\cr
0 , &\hbox{if $\a(\varepsilon)=1$.}\,\,\,\cr
\endcases
\Eq(1.theo1.1)
$$

\item{(iii)}
Let $r_n$ be an extreme space scale. For all $\b_c<\beta<\infty$, for all $\rho>0$, in $\P$-probability,
$$
\lim_{t\rightarrow 0}\lim_{n\rightarrow\infty}\CC_{n}(t,\rho t)=\asl_{\a}(1/1+\rho)\,.
\Eq(1.theo1.2)
$$
}

The next theorem shows that on extreme space (equivalently time) scales,
taking the limit $n\rightarrow \infty$ first, the process reaches
stationarity as $t\rightarrow\infty$ in the sense that the limiting
time-time correlation function is the same as that of the process $X_n$
started in its invariant measure, $\GG_n$. To state this we need a little notation.
Let $\PRM(\mu)$ be the Poisson random measure on $(0,\infty)$
with mean measure
defined through
$$
\mu(x,\infty)=x^{-\a}\,,\quad x>0\,,
\Eq(1.theo2.0)
$$
and with marks $\{\g_k\}$. Define
$$
{\CC}^{sta}_{\infty}(s)=
\sum_{k=1}^{\infty}\frac{\g_k}{\sum_{k=1}^{\infty}\g_k}e^{-s/\g_k}\,,\quad s\geq 0\,.
\Eq(1.theo2.1)
$$

\theo{\TH(1.theo2)} [Crossover to stationarity.] {\it
Let $r_n$ be an extreme space scale.
The following holds for all $\b>\b_c$.

\item{(i)} If  $\mu_n=\GG_n$ then, for all $s\leq t<t+s$,
$$
\lim_{n\rightarrow\infty}{\CC}_n(t,s)
\overset{d}\to=
{\CC}^{sta}_{\infty}(s)\,,
\Eq(1.theo2.2)
$$
where $\overset{d}\to=$ denotes equality in distribution.

\item{(ii)} If  $\mu_n=\pi_n$, for all $s\geq 0$,
$$
\lim_{t\rightarrow\infty}\lim_{n\rightarrow\infty}{\CC}_n(t,s)\overset{d}\to
={\CC}^{sta}_{\infty}(s)\,.
\Eq(1.theo2.3)
$$
}

Let us now put Theorem \thv(1.theo1) in the context of earlier results.
The very first aging results for the REM where obtained in
\cite{BBG1, BBG2}
for a discrete time version of RHT dynamics considered here, on extreme scales.
Rather than taking the double limiting procedure of \eqv(1.1.14)
we considered a decreasing sequence of extreme time scales, $c_n(E)$, $E<0$,
\note{
Specifically,
$
c_n(E)\sim e^{\b\sqrt{n}u_n(E)}
$
where $u_n(E)$ is defined through
$
2^n\P(\t_n(x)\geq c_n(E))\sim\exp\{-e^{-E}\}
$.
Thus, by \eqv(1.1') and \eqv(1.3) of Definition \thv(1.def1), $r_n\equiv c_n(E)$ is an extreme space scale.
}
and proved a statement of the form
$
\lim_{t\rightarrow\infty}\lim_{E\rightarrow -\infty}\lim_{n\rightarrow\infty}
\CC_{n, E}(t,\rho t)=\asl_{\a}(1/1+\rho)\,,
$
in $\P$-probability (see Theorem 1 of \cite{BBG2}).
The general approach of the proof is that of renewal theory; technically, it
relies on a refined knowledge of the metastable behavior of the process.
The aging scheme based on the arcsine law for subordinators was proposed later
in the landmark paper \cite{BC} and applied to the study of intermediate scales,
using potential theoretic tools.
The existence of a $\P$-almost sure arcsine aging regime was proved
in \cite{BC} (see Theorem 3.1), on the subset
$
\left\{\sfrac{3}{4}\leq\varepsilon<1 ; 0<\a(\varepsilon)<1\right\}
$
of the time scale and temperature domain, and later extended to the domain
$
\left\{0<\varepsilon<1 ; 0<\a(\varepsilon)<1\right\}
$
in \cite{CG} (see Theorem 2.1). This still leaves out the case $\varepsilon=1$
of the longest intermediate scales before extreme scales.
%
%
Let us finally stress that arcsine aging is believed to be universally present in dynamics of mean-field spin glasses.
This fact is strongly supported in \cite{BBC} where the existence of an arcisne aging regime is for
the first time proved in a model with correlations, namely the $p$-spin SK model, albeit only ``in distribution''.

So far we said nothing about the nature of the convergence of the random time-time correlation function.
Consider the last two assertions of Theorem \thv(1.theo1).
We see that statements that are valid $\P$-almost surely are turned
into in $\P$-probability statements across the line $\varepsilon=1$
(more precisely, when the relation $\frac{2^{m_n}}{2^n}\log n=o(1)$ fails to be verified).
A way to understand this transition is through the nature of the set
$
\TT_n=\{x\in\VV_n \mid \P(\t_n(x)\geq r_n)\}
$.
When $r_n$ is an intermediate space scale with $\varepsilon<1$, $\TT_n$ simply is a huge,
exponentially large site percolation cloud. In contrast, when $r_n$ is an extreme space scale, $\TT_n$,
is the very small set made of extreme $\t_n(x)$'s, and this set has a particular probabilistic structure,
namely, it resembles (for large enough $n$) a Poisson point process.
The change in the nature our convergence statements reflects this change in the nature of $\TT_n$.
This also explains why, technically, dealing with extreme scales (or very long intermediate scales)
is intrinsically more arduous then dealing with the shorter intermediate scales, below the line
$\varepsilon=1$.

To conclude our discussion of Theorem \thv(1.theo1) let us comment on the two
boundary cases $\varepsilon=0$
and $\a(\varepsilon)=1$, $0<\varepsilon\leq 1$,
where the limiting time-time correlation function
is trivial \note{
See also the paragraph after Proposition \thv(1.prop1).
}.
In both cases
one might expect that some non linear normalization of
the clock process in \eqv(1.1.8)
will produce a non trivial limit.
This was just shown to be true for short space scales where a new aging regime, coined {\it extremal aging}, had
been identified \cite{BGu, Gu}.
The critical temperature case, where a phenomenon called ``dynamical ultrametricity"
is expected to take place \cite{BB}
will be treated elsewhere.


\bigskip
\line{\bf 1.3. Convergence of the clock process.\hfill}

We now state the results on the clock process
from which Theorem \thv(1.theo1) and Theorem \thv(1.theo2) will be deduced.
We will see that to each time scale $c_n$ there correspond
{\it  auxiliary time scales} $a_n$ such that the re-scaled process
$$
S_n(t)=c_n^{-1}\wt S_n(\lfloor a_n t\rfloor)\,,\quad t\geq 0\,,
\Eq(G1.3.2')
$$
converges weakly to a subordinator.
A main aspect of our method of proof is that the limiting subordinator is constructed explicitely.
This will allow us to conclude,
comparing to the results of \cite{G1}, that both on
intermediate and extreme scales the limit is the same as in Bouchaud's
REM-like trap model
\note{
Compare Proposition \thv(1.prop1) and Proposition \thv(1.prop2) below with, respectively,
Proposition 4.8 and Proposition 4.9 of \cite{G1}.
}.
This in turn implies that on those scales the two models have
the same arcsine aging regimes.

Throughout the rest of this paper the arrow $\Rightarrow$ denotes weak convergence in the space $D([0,\infty))$
of c\`adl\`ag functions on $[0,\infty)$ equipped with the Skorohod $J_1$-topology\note{
see e\.g\. \cite{W} p\. 83 for the definition of convergence in $D([0,\infty))$.
}.
We first settle the degenerate case of short scales.
Let $\d_{\infty}$ denote the Dirac point mass at infinity.

\proposition{\TH(1.prop0)} [Short scales.] {\it
Take $\mu_n=\pi_n$ and choose $a_n$ s.t\. $a_n\sim b_n$.
Given any $0<\beta<\infty$, if $r_n$ is a short space scale that satisfies
$
\frac 1{\b  n}\log r_n\geq 4\sqrt{\frac{\log n}{n}}
$,
then, $\P$-almost surely,
$$
S_n(\cdot)\Rightarrow S^{short}(\cdot)\,,
\Eq(1.prop0.1)
$$
where $S^{short}$ is the degenerate subordinator of L\'evy measure
$
\nu^{short}=\d_{\infty}
$.
}

Turning to intermediate scales, we have:

\proposition{\TH(1.prop1)} [Intermediate scales.] {\it
Take $\mu_n=\pi_n$ and choose $a_n$ s.t\. $a_n\sim b_n$.
If $r_n$ is an intermediate space scale then, for all
$0<\varepsilon\leq 1$ and all  $0<\beta<\infty$ such that  $0<\a(\varepsilon)\leq 1$,
the following holds: $\P$-almost surely if $(2^{m_n}\log n)/2^n=o(1)$ and in $\P$-probability otherwise,
$$
S_n(\cdot)\Rightarrow S^{int}(\cdot)\,,
\Eq(1.prop1.2)
$$
where $S^{int}$ is the subordinator whose L\'evy measure,
$
\nu^{int}
$,
is defined on $(0,\infty)$ through
$$
\nu^{int}(u,\infty)=
u^{-\a(\varepsilon)}\a(\varepsilon)\G(\a(\varepsilon))\,,\quad u>0\,.
\Eq(1.prop1.1)
$$
}

Note that both $S^{short}$ and $S^{int}$ are stable subordinators of index, respectively,
$\a(0)=0$ and $0<\a(\varepsilon)\leq 1$. The cases $\a(0)=0$ and $\a(\varepsilon)= 1$ are
said to be degenerate
in the sense that the range of the subordinator reduces to the single point $0$ for the former,
and is made of the entire positive line $[0, \infty)$ for the latter.

In the final case of extreme scales the limiting
process no longer is a stable
subordinator. Neither is it a deterministic process.
As in \cite{G1}
this case is by far the more involved one\note{
See Proposition 4.9 and Section 7 of  \cite{G1}
}.
Recall that for $\mu$ defined in \eqv(1.theo2.0),
$\{\g_k\}$ denote the marks of $\PRM(\mu)$, and
introduce the re-scaled landscape variables:
$$
\g_n(x)=r_n^{-1}\t(x)\,,\quad x\in\VV_n\,.
\Eq(4.prop3.1)
$$

\proposition{\TH(1.prop2)} [Extreme scales.] {\it
Take $\mu_n=\pi_n$ and choose $a_n$ s.t\. $a_n\sim b_n$.
If $r_n$ is an extreme space scale then both the sequence of re-scaled landscapes
$(\g_n(x),\, x\in\VV_n)$,
$n\geq 1$, and the marks of $\PRM(\mu)$ can be represented on a common probability space
$(\O, \FF, \bold P)$ such that, in this representation,
denoting by $\bold S_n$ the  corresponding re-scaled clock process \eqv(G1.3.2'),
the following holds: for all $\b_c<\beta<\infty$, $\bold P$-almost surely,
$$
\bold S_n(\cdot)\Rightarrow S^{ext}(\cdot)\,,
\Eq(1.prop2.1)
$$
where
$
S^{ext}
$
is the subordinator whose L\'evy measure, $\nu^{ext}$,
is the random measure on $(0,\infty)$ defined on $(\O, \FF, \bold P)$ through
$$
\nu^{ext}(u,\infty)=\varepsilon'\sum_{k=1}^{\infty}e^{-u/\g_k}\,,\quad u>0\,,
\Eq(1.prop2.2)
$$
$\varepsilon'$ being defined in \eqv(1.3).
}

Although the limiting subordinator is not stable, the tail of the random L\'evy  measure $\nu^{ext}$
is regularly varying at $0^+$:

\lemma{\TH(1.lemma5)} {\it If $\b>\b_c$, then
$\bold P$-almost surely,
$
\nu^{ext}(u,\infty)\sim \varepsilon'u^{-\a}\a\G(\a)
$
as
$
u\rightarrow 0^+
$.
}


\bigskip
\line{\bf 1.4.  Key tools and strategy.\hfill}

The proofs of virtually all the results of the previous subsections rely on a key tool,
Theorem \thv(1.theo3), which we now state. This theorem simplifies results
from \cite{G1}, but does not specialize them to the REM dynamics. It applies to any Markov
chain $X_n$ with graph structure $(\VV_n,\EE_n)$, jump chain $J_n$, and holding time
parameters $(\l_n(x))_{x\in\VV_n}$. This will allow us to compare the
analysis of the
REM dynamics to that of  the REM-like trap model
and, based on our understanding of the latter,  to set up a strategy of proof.
Let us also recall here that the results of \cite{G1}, and hence Theorem \thv(1.theo3),
are based on a poweful result by R\. Durrett and S. Resnick \cite{DuRe},
that give conditions for partial sum processes of dependent random variables to converge to a
subordinator.

Let us further denote by $P_{\mu_n}$ the law of $J_n$ with initial distribution ${\mu_n}$, and by
$
p_n(x,y)
$
its one step transition probabilities.
Thus, in the setting of Subsection 1.1,
$p_n(x,y)$ is given by \eqv(1.9) and
$$
\l_n(x)=(\t_n(x))^{-1}\,,\quad\forall\, \,x\in \VV_n\,.
\Eq(1.6)
$$
Given sequences $c_n$ and $a_n$, the clock process $\wt S_n$ and the re-scale clock process $S_n$ are defined as
$$
\eqalign{
\wt S_n(k)&=\sum_{i=0}^{k}(\l_n(J_n(i)))^{-1}e_{n,i}\,,\quad k\in \N\,,
\cr
S_n(t)&=c_n^{-1}\wt S_n(\lfloor a_n t\rfloor)\,,\quad t\geq 0\,.
\cr
}
\Eq(1.1.6')
$$
We now formulate four conditions for the sequence $S_n$ to converge to a subordinator.
We state these conditions for a fixed realization of the random landscape,
i.e\. for fixed $\o\in\O^{\t}$,
and make this explicit by adding the superscript $\o$ to all landscape dependent quantities.

\smallskip

\noindent{\bf Condition (A1).}
There exists a $\s$-finite measure $\nu$ on $(0,\infty)$ satisfying
$\int_{(0,\infty)}(1\wedge u)\nu(du)<\infty$, such that $\nu(u,\infty)$ is continuous,
and such that, for all $t>0$ and all $u>0$,
$$
P_{\mu_n}^{\o}\left(
\left|
\sum_{j=1}^{\lfloor a_n t\rfloor}
\sum_{x\in\VV_n}p_n^{\o}(J_n^{\o}(j-1),x)
e^{-uc_n\l^{\o}_n(x)}
-t\nu^{\o}(u,\infty)
\right|
<\e
\right)=1-o(1)\,,\quad\forall\e>0\,.
\Eq(G1.A1)
$$

\noindent{\bf Condition (A2).}  For all $u>0$ and all $t>0$,
$$
P_{\mu_n}^{\o}\left(
\sum_{j=1}^{\lfloor a_n t\rfloor}\left[
\sum_{x\in\VV_n}p_n^{\o}(J_n^{\o}(j-1),x)
e^{-uc_n\l^{\o}_n(x)}
\right]^2
<\e
\right)=1-o(1)\,,\quad\forall\e>0\,.
\Eq(G1.A2)
$$

\noindent{\bf Condition (A3).}  There exists a sequence of functions $\e_n\geq 0$  satisfying
$
\displaystyle\lim_{\d\rightarrow 0}\limsup_{n\rightarrow \infty}\e_n(\d)=0
$
such that, for some $0<\d_0\leq 1$, for all $0<\d\leq\d_0$ and all $t>0$,
$$
E_{\mu_n}^{\o}\left(\int_{0}^{\d}du
\sum_{j=1}^{\lfloor a_n t\rfloor}\sum_{x\in\VV_n}p_n^{\o}(J_n^{\o}(j-1),x)e^{-uc_n\l^{\o}_n(x)}
\right)
\leq t\e_n(\d)\,.
\Eq(G1.A3)
$$

\noindent{\bf Condition (A0').}
$$
\sum_{x\in\VV_n}\mu_n^{\o}(x)e^{-vc_n\l^{\o}_n(x)}=o(1)\,.
\Eq(1.A0')
$$

\theo{\TH(1.theo3)}
{\it For all sequences of initial distributions  $\mu_n$ and all sequences  $a_n$ and $c_n$
for which Conditions (A0'), (A1), (A2) and (A3) are verified,
either $\P$-almost surely or in $\P$-probability,
the following holds w.r.t\. the same convergence mode:
Let $\{(t_k,\xi_k)\}$ be the points of a Poisson random measure of intensity measure
$dt\times d\nu$.  We have,
$$
 S_n(\cdot)\Rightarrow  S(\cdot)=\sum_{t_k\leq \cdot}\xi_k\,.
\Eq(G1.3.theo1.1)
$$
Moreover,
$$
\lim_{n\rightarrow\infty}\CC_{n}(t,s)=\CC_{\infty}(t,s)\,,
\Eq(G1.3.theo1.4)
$$
where
$$
\CC_{\infty}(t,s)=\PP\left(\left\{ S(u)\,,u>0\right\}
\cap (t, t+s)=\emptyset\right)\,,\quad 0\leq t<t+s\,.
\Eq(G1.3.theo1.3)
$$
}

\proof
Theorem \thv(1.theo3) results from the concatenation of Theorem 1.3 of \cite{G1}
and a trimmed version of Theorem 1.4 of \cite{G1}, where Condition (A0) is
replaced by the more restrictive Condition (A0').
\endproof

\remark Unlike in \cite{G1} we did not separate the first steps of the clock process from
the remaining ones. Instead, we made the simplifying assumption (see Condition (A0'))
that the first step converges to zero. The question of more general initial
distributions will be discussed elsewhere.

We see from Theorem \thv(1.theo3) that the L\'evy measure $\nu$
of the limiting subordinator in \eqv(G1.3.theo1.1) is determined by Condition (A1).
In order to prove that
the limiting re-scaled clock process in the REM is the same as in the REM-like trap model,
we first have to understand why and how Condition (A1) can be satisfied with the same
measure $\nu$. To this hand take $\mu_n=\pi_n$ for both models, $\pi_n$ being the invariant measure of $J_n$,
and set
$$
\nu_n^{J,t}(u,\infty)=\sum_{j=1}^{\lfloor a_n t\rfloor}\sum_{x\in\VV_n}p_n^{\o}(J_n^{\o}(j-1),x)e^{-uc_n\l^{\o}_n(x)}\,,\quad u>0\,.
\Eq(1.4.1)
$$
Now note that, by reversibility, for both models,
$$
E_{\pi_n}\nu_n^{J,t}(u,\infty)=\lfloor a_n t\rfloor \sum_{x\in\VV_n}\pi_n(x)e^{-uc_n\l^{\o}_n(x)}\,,\quad u>0\,.
\Eq(1.4.2)
$$
If we now specialize to the REM-like trap model
the following two features will suffice to determine the measure $\nu$:
\item{(i)} $(\l_n^{-1}(x), x\in\VV_n)$ are i.i.d\. r.v.'s in the domain of attraction of a positive stable law with index $0<\a<1$, and

\item{(ii)} $\nu_n^{J,t}(u,\infty)=E_{\pi}\nu_n^{J,t}(u,\infty)$.

\noindent
Indeed, from \eqv(1.4.2) and (ii) one expects that, by independence and a concentration argument,
$\nu(u,\infty)-\E E_{\pi}\nu_n^{J,t}(u,\infty)\downarrow 0$ as $n\uparrow\infty$ with
$\P$-probability approaching one;
using (i) then allows to compute $\E E_{\pi}\nu_n^{J,t}(u,\infty)$, and so
$\nu=\lim_{n\rightarrow\infty}\E E_{\pi}\nu_n^{J,t}(u,\infty)$ provided that the limit exists\note{
This reasoning applies on intermediate scales. A more refined argument is be needed on extreme scales.}.

The features (i) $\&$ (ii), that are put in by hand in REM-like trap model, are clearly not present in the
RHT dynamics of the REM. However, if we can prove that (i')  the re-scaled landscape $(c_n\l_n^{-1}(x), x\in\VV_n)$  of the REM
is heavy tailed (i.e\., show that there exists sequences $a_n$ and $c_n$ such that
$
a_n\P(\l_n^{-1}(x)>u c_n)\sim u^{-\a}
$,
for some $0<\a<1$), and  that, loosely speaking,
(ii') the quantity $\nu_n^{J,t}(u,\infty)$ obeys an ergodic theorem with $\P$-probability close to one,
then we should be able to reduce the analysis of the quantity
$\nu_n^{J,t}(u,\infty)$ of the REM to that of the REM-like trap model.
An analogous reasoning applies to Condition (A2).
If (i') is known to hold from earlier works
(at least on intermediate scales), how to justify (ii') is less immediate.
To do this we will heavily draw on the specific properties of the jump chain $J_n$,
in particular,
on the fact that its trajectories do not depend on the randomness of the lansdcape.
This is a key feature of the random hopping time dynamics dynamics, and the reason why its analysis is
so much simpler than that of the usual metropolis dynamics \cite{G2}.

The rest of this paper is organized as follows.
The heavy tailed nature of the trapping landscape is established in Section 2.
Section 3 collects the properties of the jump chain $J_n$ that will be needed in Section 4
to establish an `ergodic theorem' for each of the
sums appearing in Condition (A1) and Condition (A2)
(called respectively $\nu_n^{J,t}(u,\infty)$ and $(\s_n^{J,t})^2(u,\infty)$).
Once this is done it remains to establish the properties of the
chain averaged sums:
for this we separate the case of short and intermediate scales,
treated in Section 5, from the case of extreme scales, dealt with in Section 6.
The proofs of results obtained on short and intermediate scales (respec\., extreme scales)
are stated in Section 5, (respec\., Section 6).

\bigskip

\chap{2. Properties of the landscape.}2

In this section we establish the needed properties of the re-scaled landscape variables
$(r_n^{-1}\t_n(x), x\in\VV_n)$,
and most importantly, the  heavy tailed nature of their distribution.
The notations of Subsection 1.2 (from \eqv(1.1') to \eqv(1.theo1.0')) are in force throughout this section.
We assume that $0<\b<\infty$ is fixed, and as before, drop all dependence on $\b$ in the notation.
For $u\geq 0$ set $G_n(u)=\P(\tau_n(x)>u)$.
Since this is a continuous monotone decreasing function, it has a well defined inverse
$
G_n^{-1}(u):=\inf\{y\geq 0 : G_n(y)\leq u\}
$.

For $ v\geq 0$ set
$$
h_n(v)=b_nG_n(r_nv)\,.
\Eq(2.13)
$$

\lemma{\TH(2.lemma6)} {\it
Let $r_n$ be any of the space scales of Definition \thv(1.def1).

\item{(i)} For each fixed $\zeta>0$ and all $n$ sufficiently large so that $\zeta>r_n^{-1}$,
the following holds:
for all $v$ such that $\zeta\leq v<\infty$,
$$
h_n(v)= v^{-\a_n}(1+o(1))\,,
\Eq(2.lem6.1)
$$
where $0\leq \a_n=\a(\varepsilon)+o(1)$.

\item{(ii)} Let $0<\d<1$.
Then, for all $v$ such that $r_n^{-\d}\leq v\leq 1$ and all large enough $n$,
$$
v^{-\a_n}(1+o(1))
\leq h_n(v)\leq
\sfrac{1}{1-\d}v^{-\a_n(1-\frac{\d}{2})}(1+o(1))\,,
\Eq(6.lem11.1)
$$
where $\a_n$ is as before.
}

Next, for $u\geq 0$ set
$$
g_n(u)=r_n^{-1}G_n^{-1}(u/b_n)\,.
\Eq(2.12)
$$
It is plain that $g_n(v)=h_n^{-1}(v)$. It is plain also that both $g_n$ and $h_n$
are continuous monotone decreasing functions.
The following lemma is tailored to deal with the case of extreme space scales. Recall that $\a\equiv\a(1)$.

\lemma{\TH(2.lemma4)} {\it Let $r_n$ be an extreme space scale.

\item{(i)} For each fixed $u>0$, for any sequence $u_n$ such that
$|u_n-u|\rightarrow 0$  as $n\rightarrow\infty$,
$$
g_n(u_n)\rightarrow u^{-(1/\a)}\,,\quad n\rightarrow\infty\,.
\Eq(2.lem4.1)
$$
\item{(ii)} There exists a constant $0<C<\infty$
such that, for all $n$ large enough,
$$
g_n(u)\leq u^{-1/\a}C\,,\quad 1\leq u\leq b_n(1-\Phi(1/(\b\sqrt n)))\,,
\Eq(2.lem4.2)
$$
where $\Phi$ denotes the standard Gaussian distribution function .
}

The proofs of these two lemmata rely on Lemma \thv(2.lemma3) below.
%
%
Denote by $\Phi$ and $\phi$ the standard Gaussian distribution function and density, respectively.
Let $B_n$ be defined through
$$
b_n\frac{\phi(B_n)}{B_n}=1\,,
\Eq(2.7)
$$
and set
$
A_n=B_n^{-1}
$

\lemma{\TH(2.lemma3)} {\it Let $r_n$ be any space scale.
Let $\wt B_n$ be a sequence such that, as $n\rightarrow\infty,$
$$
\eqalign{
\d_n&:=(\wt B_n-B_n)/A_n\rightarrow 0\,.
}
\Eq(2.9)
$$
Then, for all $x$ such that $A_n x+ \wt B_n>0$ for large enough $n$,
$$
b_n(1-\Phi(A_n x+ \wt B_n))=\frac{\exp\left(-x\left[1+\sfrac{1}{2}A_n^2x\right]\right)}
{1+A_n^2x}
\left\{1+\OO\left(\d_n[1+A_n^2+A_n^2x]\right)+\OO(A_n^2)\right\}\,.
\Eq(2.10)
$$
}
\proof The lemma is a direct consequence of the following expressions, valid for all $x>0$,
$$
\eqalign{
1-\Phi(x)&=x^{-1}\phi(x)-r(x)\cr
&=x^{-1}(1-x^{-2})\phi(x)-s(x)\,,
}
\Eq(2.11)
$$
where
$$
\eqalign{
0<r(x)=\int_x^{\infty}y^{-2}\phi(y)dy< x^{-3}\phi(x)\,,\cr
0<s(x)\int_x^{\infty}y^{-4}\phi(y)dy< x^{-5}\phi(x)\,,
}
\Eq(2.11')
$$
(see \cite{AS}, p\. 932). We leave the details to the reader.
\endproof

We now prove Lemmata \thv(2.lemma6) and \thv(2.lemma4), beginning with Lemma \thv(2.lemma6).

\proofof{Lemma \thv(2.lemma6)}
By definition of $G_n$ we may write
$$
h_n(v)
=b_n\Bigl(1-\Phi\left(A_n\log(v^{\a_n})+\overline B_n\right)\Bigr)\,,
\Eq(2.14)
$$
where
$$
\eqalign{
\overline B_n&=(\b\sqrt n)^{-1} \log r_n\,,\cr
\a_n&=(\b\sqrt n)^{-1}B_n\,.
}
\Eq(2.15)
$$

We first claim that for $v\geq r_n^{-1}$, which guarantees that $A_n\log(v^{\a_n})+\overline B_n>0$,  the  sequence $\overline B_n$
satisfies the assumptions
of Lemma \thv(2.lemma3).
For this we use the know fact that the sequence $\wh B_n$ defined by
$$
\wh B_n=(2\log b_n)^{\frac{1}{2}}-\sfrac{1}{2}(\log\log b_n +\log 4\pi)/(2\log b_n)^{\frac{1}{2}}\,,
\Eq(2.18)
$$
satisfies
$$
(\wh B_n-B_n)/A_n=o(1)
\Eq(2.19)
$$
(see \cite{H}, p\. 434, paragraph containing Eq\. (4)).
By \eqv(2.11)-\eqv(2.11') we easily get that
$$
b_n\bigl(1-\Phi(\wh B_n)\bigr)=1-(\log\log b_n)^2(16\log b_n)^{-1}(1+o(1))\,,
\Eq(2.20)
$$
whereas, by definition of $b_n$ (see \eqv(1.1')),
$$
b_n\bigl(1-\Phi(\overline B_n)\bigr)=1\,.
\Eq(2.21)
$$
Since $\Phi$ is monotone and increasing, \eqv(2.20) and \eqv(2.21) imply that $\wh B_n>\overline B_n$.
Thus
$$
\bigl(1-\Phi(\overline B_n)\bigr)-\bigl(1-\Phi(\wh B_n)\bigr)
=\Phi(\wh B_n)-\Phi(\overline B_n)
\geq\phi(\wh B_n)(\wh B_n-\overline B_n)\geq 0\,.
\Eq(2.22)
$$
This, together with \eqv(2.20) and \eqv(2.21), yields
$$
\eqalign{
0<\wh B_n-\overline B_n
&<\bigl[b_n\phi(\wh B_n)\bigr]^{-1}(\log\log b_n)^2(16\log b_n)^{-1}(1+o(1))\,.
\cr
}
\Eq(2.23)
$$
Now, by \eqv(2.7),
$$
b_n\phi(\wh B_n)=B_n\bigl[\phi(\wh B_n)/\phi(B_n)\bigr]
=B_n\exp\bigl\{-\sfrac{1}{2}(\wh B_n-B_n)(\wh B_n+ B_n)\bigr\}
\leq B_n(1+o(1))\,,
\Eq(2.24)
$$
where the final inequality follows from \eqv(2.19).
Finally, combining \eqv(2.23) and \eqv(2.24) yields
$
0<(\wh B_n-\overline B_n)/A_n=o(1)
$,
and using again \eqv(2.19),
we obtain
$
(\overline B_n-B_n)/A_n=o(1)
$,
which was the claim.

To control the behavior of the sequences $A_n$, $\a_n$, and $r_n$, we will need an expression for the
(of course well known) solution $B_n$ of \eqv(2.7). Here is one (\cite{Cr}, p\. 374):
$$
B_n=(2\log b_n)^{\frac{1}{2}}-\sfrac{1}{2}(\log\log b_n +\log 4\pi)/(2\log b_n)^{\frac{1}{2}}+\OO(1/\log b_n)\,.
\Eq(2.16)
$$
Note that so far we didn't make use of the assumption on $r_n$. The choice of
$r_n$ steps in only now: using \eqv(2.15), \eqv(2.16), the fact that
$
2\log b_n=(2\log2)m_n=\b_c^2(m_n/n)n
$,
and the just established fact that
$
(\overline B_n-B_n)/A_n\rightarrow 0
$, we obtain,
for intermediate space scales,
$$
\eqalign{
\log r_n&=\b\b_c(\varepsilon) n(1+o(1))\,,\cr
\log b_n&=\sfrac{1}{2}\b^2_c(\varepsilon)n(1+o(1))\,,\cr
\a_n&=\a(\varepsilon)(1+o(1))\,,
}
\Eq(2.23'')
$$
and, for short space scales,
$$
\eqalign{
\log r_n&=o(n)\,,\cr
\log b_n&=o(n)\,,\cr
\a_n&=o(1)\,.
}
\Eq(2.23''')
$$
Finally for extreme space scales, wrting $\b_c(1)\equiv\b_c$, we have that
$
2\log b_n=(2\log2)m_n=\b_c^2n(1-C/n)
$
for some constant $0<C<\infty$. Thus, instead of \eqv(2.23''), we get:
$$
\eqalign{
\log b_n&=\sfrac{1}{2}\b^2_cn(1-C/n)\,,\cr
\log r_n&=\b\b_c n(1-o(1))\,,\cr
\a_n\leq \a&\text{and}
\a_n=\a(1-o(1))\,.
}
\Eq(2.23')
$$

We are now equipped to prove Lemma \thv(2.lemma6).
By Lemma \thv(2.lemma3), for all  $v\geq r_n^{-1}$,
$$
h_n(v)=\frac{\exp\left(-\a_n\log v\left[1+\sfrac{1}{2}A_n^2\a_n\log v\right]\right)}
{1+A_n^2\a_n\log v}
\left\{1+\OO\left(\d_n[1+A_n^2+A_n^2\a_n\log v]\right)+\OO(A_n^2)\right\}\,,
\Eq(2.24')
$$
where $\d_n\downarrow 0$ as $n\uparrow\infty$. Therefore, for each fixed $0< v<\infty$,
and all large enough $n$ so that $v> r_n^{-1}$,
$$
h_n(v)=v^{-\a_n}(1+o(1))\,.
\Eq(2.25bis)
$$
This together with \eqv(2.23'') and \eqv(2.23''') proves assertion (i) of the lemma.
To prove  assertion (ii) note that by \eqv(2.15), since $A_n=B_n^{-1}$,
$
A_n^2\a_n=\frac{1}{\log r_n}\frac{\overline B_n}{B_n}
$
where $\frac{\overline B_n}{B_n}=1+o(1)$ (see the paragraph following \eqv(2.24)).
Thus,
for all $v$ satisfying $r_n^{-\d}\leq v\leq 1$, we have
$$
-\d\sfrac{\overline B_n}{B_n}\leq A_n^2\a_n\log v\leq 0\,.
\Eq(2.25ter)
$$
Combining this and \eqv(2.24') immediately yields the bounds \eqv(6.lem11.1).
The proof of Lemma \thv(2.lemma6) is now done.
\endproof


\proofof{Lemma \thv(2.lemma4)}
Up until \eqv(2.24') we proceed exactly as in the proof of Lemma \thv(2.lemma6).
Now, by \eqv(2.24'), for each fixed $0\leq v<\infty$,
any sequence $v_n$ such that $|v_n-v|\rightarrow 0$, and all large enough $n$ (so that $v> r_n^{-1}$),
$$
h_n(v_n)=v_n^{-\a_n}(1+o(1))= v^{-\a(1-o(1))}(1+o(1))\,.
\Eq(2.25)
$$
This and the relation $g_n(v)=h_n^{-1}(v)$  imply that for each fixed $0<u<\infty$,
any sequence $u_n$ such that $|u_n-u|\rightarrow 0$, and all large enough $n$ (so that $u<h_n(r_n^{-1})$),
$$
g_n(u_n)=u_n^{-(1/\a_n)}(1+o(1))= u^{-(1/\a)(1+o(1))}(1+o(1))\,,
\Eq(2.26)
$$
which is tantamount to assertion (i) of the lemma.

To prove assertion (ii) assume that $r_n^{-1}\leq v\leq 1$.
Recall that $h_n$ is a monotonous function so that
if $h_n(v)=g_n^{-1}(v)$ for all $r_n^{-1}\leq v\leq 1$,
then $g_n(u)=h_n^{-1}(u)$ for all $h_n(1)\leq u\leq h_n(r_n^{-1})$.
Now $h_n(1)=b_nG_n(r_n)=1$, as follows from \eqv(1.1'), and
$h_n(r_n^{-1})=b_nG_n(1)=b_n(1-\Phi(1/(\b\sqrt n)))$.
Observe next that
$r_n^{-1}\leq v\leq 1$ is equivalent to
$
-1\leq A_n^2\log v^{\a_n} \leq 0
$.
Therefore, by \eqv(2.24'), for large enough $n$,
$$
h_n(v)\geq (1-2\d_n)v^{-\a_n}\,,\quad r_n^{-1}\leq v\leq 1\,.
\Eq(2.27)
$$
By monotonicity of $h_n$,
$$
g_n(u)=h_n^{-1}(u)\leq (1-2\d_n)^{1/\a_n} u^{-1/\a_n}\,,\quad 1\leq u\leq b_n(1-\Phi(1/(\b\sqrt n)))\,.
\Eq(2.28)
$$
From this and the fact that
$\a_n\leq \a$ (see \eqv(2.23')), \eqv(2.lem4.2) readily obtains.
This concludes the proof of the lemma.
\endproof

\remark We see from the proof of Lemma \thv(2.lemma4) that the lemma holds true not only for extreme scales,
but for intermediate scales also provided one replaces $\a$ by $\a(\varepsilon)$ everywhere.

\bigskip


\chap{3. The jump chain: some estimates.}3

This section is about the jump chain $J_n$, i.e\. the simple random walk. We gather here all the results
that will be needed later
to reduce Condition (A1) and Condition (A2)
of Theorem \thv(1.theo3) to  conditions that are independent from $J_n$.
Proposition \thv(3.prop1) below is central to this scheme.
It will allow us to substitute
the stationary chain for the jump chain after $\theta_n\sim n^2$ steps have been taken.

\proposition{\TH(3.prop1)} {\it Set
$
\theta_n=2
\left\lceil
\frac{3}{2}(n-1)\log 2/\left|\log\left(1-\frac{2}{n}\right)\right|
\right\rceil
$.
For all pairs $x\in\VV_n, y\in\VV_n$, and all $i\geq 0$,
$$
P_{\pi_n}\left(J_n(i+\theta_n)=y,  J_n(0)=x \right)+P_{\pi_n}\left(J_n(i+1+\theta_n)=y, J_n(0)=x \right)
=2(1+\d_n)\pi_n(x)\pi_n(y)
\,,
\Eq(3.prop.1)
$$
where $|\d_n|\leq 2^{-n}$.
}

\noindent The next two propositions are technical estimates
needed in the proofs of
Proposition \thv(4.prop1) and Proposition \thv(6.prop1) respectively.
Let $p_n^{l}(\cdot,\cdot)$ denote the $l$ steps transition probabilities of $J_n$,
and let $\dist(x,y):=\#\bigl\{i\in\{1,\dots,n\}\,:\,x_i\neq y_i\bigr\}$
be the Hamming distance.

\proposition{\TH(3.prop3)} {\it For all $m\leq n^2$,
$$
\sum_{l=1}^{2m}p_n^{l+2}(z,z)\leq \frac{c}{n^2}\,,\quad\forall z\in\VV_n\,,
\Eq(3.prop3.1)
$$
for some constant $0<c<\infty$.
}

\proposition{\TH(3.prop4)} {\it For all $m\leq n^2$, for all pairs of distinct vertices $y,z\in\VV_n$
satisfying $\dist(y,z)=\frac{n}{2}(1-o(1))$,
$$
\sum_{l=1}^{2m}p_n^{l+2}(y,z)\leq e^{-cn}\,,
\Eq(3.prop4.1)
$$
for some constant $0<c<\infty$.
}

We first prove Proposition \thv(3.prop1).
For this we will use the following classical bound by Diaconis and Stroock \cite{DS}
for the total variation distance to stationarity of reversible Markov chains.
Recall that the total variation distance between two probabilities $\mu$ and $\mu'$ on $\VV^{\circ}_n$ is defined as
$$
\|\mu-\mu'\|_{TV}=\max_{A\subset\VV^{\circ}_n} |\mu(A)-\mu'(A)|=\frac{1}{2}\sum_{x\in\VV^{\circ}_n}|\mu(x)-\mu'(x)|\,.
$$

\proposition{\TH(3.prop2)} [\cite{DS}, Proposition 3] {\it Let $\nu$, $Q(x,y)$
be a reversible Markov chain on a finite set $X$. Assume that the one-step transition
probability matrix $Q$ is irreducible with eigenvalues $1=\b_0>\b_1\geq\dots \geq\b_{k-1}\geq -1$.
Then, for all $x\in X$ and $m\in\N$,
$$
4\bigl\|Q^m(x,\cdot)-\nu(\cdot)\bigr\|^2_{TV}\leq \sfrac{1-\nu(x)}{\nu(x)}\b_*^{2m}\,,
\quad\b_*=\min(\b_1,|\b_{k-1}|)\,.
\Eq(3.prop2.1)
$$
}

Proposition \thv(3.prop2) cannot be applied directly
to the jump chain $J_n$, since it is periodic with period two, but
it can be applied to the aperiodic chains obtained by observing $J_n$ at even, respectively, odd times.
In view of doing this let us partition the cube into the sub-cubes
$
\VV_n=\VV_n^{od}\cup\VV_n^{ev}
$
of vertices that are at odd, respectively, even distance of the vertex $x=(1,1\dots,1)$:
$$
\eqalign{
\VV_n^{od}&=\{x\in\VV_n \mid \textstyle{ \sum_{i=1}^n} (x_i+1)/2 \,\,\, \hbox{\rm is odd}\,\}\,,\cr
\VV_n^{ev}&=\{x\in\VV_n \mid \textstyle{ \sum_{i=1}^n} (x_i+1)/2\,\,\, \hbox{\rm is even}\,\}\,.\cr
}
\Eq(3.5)
$$
To each of these sub-cubes we associate a chain, $J_n^{od}$ and $J_n^{ev}$,
as follows.
Let the symbol ${}^{\circ}$ denote either ${\it od}$ or ${\it ev}$.
Set $\EE_n^{\circ}=\left\{(x,y)\in\VV^{\circ}_n\times\VV^{\circ}_n\,:\, \sum_{i=1}^n|x_i-y_i|=4 \right\}$.
Then $(J_n^{\circ}(k)\,,k\in\N)$  is the chain on $\VV^{\circ}_n$ with transition probabilities
$
p_n^{\circ}(x,y)=P^{\circ}(J^{\circ}_n(i+1)=y\mid J^{\circ}_n(i)=x)=P_{\pi_n}(J_n(i+2)=y\mid J_n(i)=x)
$,
that is,
$$
p_n^{\circ}(x,y)=
\cases
\frac{2}{n^2}&\hbox{if $(x,y)\in\EE_n^{\circ}$} , \,\,\,\cr
\frac{1}{n}&\hbox{if $x=y$} , \,\,\,\cr
0 , &\hbox{otherwise.}\,\,\,\cr
\endcases
\Eq(3.6)
$$
Clearly $J_n^{\circ}$ is aperiodic and has a unique reversible invariant measure $\pi_n^{\circ}$ given by
$$
\pi_n^{\circ}(x)=
\cases
2^{-n+1}&\hbox{if $x\in\VV_n^{\circ}$} , \,\,\,\cr
0 , &\hbox{otherwise.}\,\,\,\cr
\endcases
\Eq(3.7)
$$
In what follows we denote by $P^{\circ}$ the law of $J^{\circ}_n$ with initial distribution $\pi^{\circ}_n$.

Applying Proposition \thv(3.prop2) to each of the two chains $J^{\circ}_n$ yields:

\lemma{\TH(3.lemma1)} {\it Let $\theta_n$ and $\d_n$ be as in Proposition \thv(3.prop1).
Then, for all $x\in\VV^{\circ}_n$, all $y\in\VV_n$, and  large enough  $n$,
$
P^{\circ}\left(J^{\circ}_n(l)=y\mid J^{\circ}_n(0)=x\right)=(1+\d_n)\pi^{\circ}_n(y)
$,
for all  $l\geq\theta_n/2$.
}

\proofof{Lemma \thv(3.lemma1)}
The eigenvalues  of the transition matrix
$Q^{\circ}:=\left(p^{\circ}_n(x,y)\right)_{\VV^{\circ}_n\times\VV^{\circ}_n}$ of $J^{\circ}_n$
are $\left(1-2\frac{j}{n}\right)^2$, $0\leq j\leq \lfloor \frac{n}{2}\rfloor$,
where ${}^{\circ}$ denotes either ${\it od}$ or ${\it ev}$. The proof of this statement uses
the following three facts: (i) firstly, the eigenvalues of the transition matrix
$Q:=\left(p_n(x,y)\right)_{\VV_n\times\VV_n}$ of $J_n$ are $1-2j/n$, $0\leq j\leq n$
(see, for example, \cite{DS} example 2.2 p\. 45);
(ii) next,  by definition of $Q^{\circ}$, $Q^2=Q^{ev}+Q^{od}$ and  $Q^{ev}Q^{od}=Q^{od}Q^{ev}=0$;
(iii) finally, $Q^{ev}$ and $Q^{od}$ can be obtained from one another by permutation of their
rows and columns. Now it follows from (iii) that $Q^{ev}$ and $Q^{od}$ must have the same spectrum.
This fact combined with (i) and (ii) imply that this spectrum coincide with that of $Q^2$.
The conclusion then follows from (i).

We may thus apply \eqv(3.prop2.1) to the chain  $J^{\circ}_n$ with $\b_*=\left(1-\frac{2}{n}\right)^2$.
Choosing $m=\theta_n/2
=\left\lceil
\frac{3}{2}(n-1)\log 2/\left|\log\left(1-\frac{2}{n}\right)\right|
\right\rceil
$, this yields
$
P^{\circ}\left(J^{\circ}_n(l)=y\mid J^{\circ}_n(0)=x\right)=(1+\d_n)\pi^{\circ}_n(y)
$
where $\d_n^2\leq \frac{1}{4}2^{3(n-1)}\left(1-\frac{2}{n}\right)^{2\theta_n}\leq 2^{-3n+1}$
for all $n$ large enough, and thus $|\d_n|\leq  2^{-n}$.
The lemma is proven.\endproof

\proofof{Proposition \thv(3.prop1)}
Proposition \thv(3.prop1) is an immediate consequence of Lemma \thv(3.lemma1).
To see this set
(to simplify the notation we drop the dependence of $P$ on its initial distribution $\pi_n$ in this proof)
$
\D=P\left(J_n(i+\theta_n)=y,  J_n(0)=x \right)+P\left(J_n(i+1+\theta_n)=y, J_n(0)=x \right)
$.
Since $J_n$ is started from its invariant measure $\pi_n$,
$
\D=P\left(J_n(i+\theta_n)=y\mid J_n(0)=x \right)\pi_n(x)+P\left(J_n(i+1+\theta_n)=y\mid J_n(0)=x \right)\pi_n(x)
$.
Without loss of generality we may assume that $i+\theta_n$ is even, and set $i+\theta_n=2l$.
Then, using the notation $x\sim y\Leftrightarrow\dist(x,y)=1$, $\D/\pi_n(x)$ can we rewritten as
$$
\eqalign{
\frac{\D}{\pi_n(x)}
&=P\left(J_n(2l)=y\mid J_n(0)=x \right)+P\left(J_n(2l+1)=y\mid J_n(0)=x \right)
\cr
=&P\left(J_n(2l)=y\mid J_n(0)=x \right)+\frac{1}{n}\sum_{z\sim x}P\left(J_n(2l+1)=y\mid J_n(1)=z \right)
\cr
=&
P\left(J_n^{ev}(l)=y\mid J_n^{ev}(0)=x \right)\1_{\{x\in\VV_n^{ev}\}}
+P\left(J_n^{od}(l)=y\mid J_n^{od}(0)=x \right)\1_{\{x\in\VV_n^{od}\}}
\cr
+&\frac{1}{n}\sum_{z\sim x}\Bigl[
P\left(J_n^{od}(l)=y\mid J_n^{od}(1)=z \right)\1_{\{x\in\VV_n^{ev}\}}
+P\left(J_n^{ev}(l)=y\mid J_n^{ev}(1)=z \right)\1_{\{x\in\VV_n^{od}\}}
\Bigr]
\,,
}
\Eq(3.prop.3)
$$
and, making use of Lemma \thv(3.lemma1),
$$
\eqalign{
\frac{\D}{\pi_n(x)}
=(1+\d_n)\Bigl[
\pi_n^{ev}(y)\1_{\{x\in\VV_n^{ev}\}}+\pi_n^{od}(y)\1_{\{x\in\VV_n^{od}\}}
+\pi_n^{od}(y)\1_{\{x\in\VV_n^{ev}\}}+\pi_n^{ev}(y)\1_{\{x\in\VV_n^{od}\}}
\Bigr]\,.
}
\Eq(3.prop.4)
$$
Now, exactly one of the indicator function in the right hand side of \eqv(3.prop.3) is non zero, so that,
by \eqv(3.7) and \eqv(1.8),
$
\D=2(1+\d_n)\pi_n(x)\pi_n(y)
$.
The proof of Proposition \thv(3.prop1) is done.\endproof

We now prove Proposition \thv(3.prop3) and Proposition \thv(3.prop4).

\proofof{Proposition \thv(3.prop3)} Consider the Ehrenfest chain on state space $\{0,\dots,2n\}$
with one step transition probabilities $r_n(i,i+1)=\frac{i}{2n}$ and $r_n(i, i-1)=1-\frac{i}{2n}$.
Denote by $r_n^{l}(\cdot,\cdot)$ its $l$ steps transition probabilities.
It is well known (see e.g. \cite{BG}) that
$
p_n^{l}(z,z)=r_n^{l}(0,0)
$
for all $l\geq 0$ and all $z\in\VV_n$.
Hence
$
\sum_{l=1}^{2m}p_n^{l+2}(z,z)=\sum_{l=1}^{2m}r_n^{l+2}(0,0)
$.
It is in turn well known (see \cite{Kem}, p\. 25, formula (4.18)) that
$$
r_n^{l}(0,0)=2^{-n}\sum_{k=0}^{2n}\binom{2n}{l}\left(1-\frac{k}{n}\right)^l\,,\quad l\geq 1\,.
\Eq(3.prop3.2)
$$
Note that by symmetry, $r_n^{2l+1}(0,0)=0$. Simple calculations yield
$r_n^{4}(0,0)=\frac{c_2}{n^2}$, $r_n^{6}(0,0)=\frac{c_3}{n^3}$, and $r_n^{8}(0,0)=\frac{c_4}{n^4}$,
for some constants $0<c_i<\infty$, $1\leq i\leq 3$. Thus, if $m\leq 3$,
$
\sum_{l=1}^{2m}r_n^{l+2}(0,0)\leq \frac{c}{n^2}
$
for some constant $0<c<\infty$. If now $m>3$, write
$
\sum_{l=1}^{2m}r_n^{l+2}(0,0)
=r_n^{4}(0,0)+r_n^{6}(0,0)+\sum_{l=6}^{2m}r_n^{l+2}(0,0)
$,
and use that by \eqv(3.prop3.2),
$$
\sum_{l=6}^{2m}r_n^{l+2}(0,0)
=
2^{-n}\sum_{k=0}^{2n}\binom{2n}{l}\sum_{l=6}^{2m}\left(1-\sfrac{k}{n}\right)^{l+2}
\leq
2^{-n}\sum_{k=0}^{2n}\binom{2n}{l}\left(1-\sfrac{k}{n}\right)^{8}\sum_{j=0}^{m-1}\left(1-\sfrac{k}{n}\right)^j\,.
\Eq(3.prop3.3)
$$
Since $|1-\sfrac{k}{n}|\leq 1$,
$
\sum_{l=6}^{2m}r_n^{l+2}(0,0)
\leq 2^{-n}\sum_{k=0}^{2n}\binom{2n}{l}\left(1-\sfrac{k}{n}\right)^{8}m
=mr_n^{8}(0,0)
\leq n^2\sfrac{c_4}{n^4}
$,
so that
$
\sum_{l=1}^{2m}r_n^{l+2}(0,0)
\leq \frac{c_2}{n^2}+\frac{c_3}{n^3}+n^2\sfrac{c_4}{n^4}
\leq \sfrac{c}{n^2}
$
for some constant $0<c<\infty$. The lemma is proven.\endproof

\proofof{Proposition \thv(3.prop4)}
This estimate is proved using  a $d$-dimensional version of the Ehrenfest scheme known as
the lumping procedure, and studied e.g\. in \cite{BG}.
In what follows we mostly stick to the notations of \cite{BG}, hoping that this will create
no confusion. Without loss of generality we may take $y\equiv 1$
to be the vertex whose coordinates all take the value 1. Let $\g^{\L}$ be the map
(1.7) of \cite{BG} derived from the partition of $\L\equiv\{1,\dots,n\}$ into
$d=2$ classes, $\L=\L_1\cup\L_2$, defined through the relation:
$i\in\L_1$ if the $i^{th}$ coordinate of $z$ is 1, and $i\in\L_2$ otherwise.
 The resulting lumped chain $X^{\L}_n\equiv\g^{\L}(J_n)$
has range $\G_{n,2}=\g^{\L}(\VV_n)\subset[-1,1]^2$. Note that the vertices
1 and $y$ of $\VV_n$ are mapped respectively on the corners $1\equiv(1,1)$ and $x\equiv(1,-1)$ of $[-1,1]^2$.
Without loss of generality we may assume that $0\in\G_{n,2}$.
Now, denoting by $\P^{\circ}$ the law of $X^{\L}_n$, we have,
$
p_n^{l+2}(y,z)=\P^{\circ}(X^{\L}_n(l+2)=x \mid X^{\L}_n(0)=1)
$.
Let
$
\t^{x'}_x=\inf\{ k>0 \mid X^{\L}_n(0)=x', X^{\L}_n(k)=x\}
$.
Starting from 1, the lumped chain may visit 0 before it visits x or not. Thus
$
p_n^{l+2}(1,z)
=\P^{\circ}(X^{\L}_n(l+2)=x, \t^{1}_0< \t^{1}_x)+\P^{\circ}(X^{\L}_n(l+2)=x, \t^{1}_0\geq \t^{1}_x)
$.
On the one hand, using Theorem 3.2 of \cite{BG},
$
\P^{\circ}(X^{\L}_n(l+2)=x, \t^{1}_0\geq \t^{1}_x)
\leq \P^{\circ}(\t^{1}_x\leq \t^{1}_0)
\leq e^{-c_1n}
$
for some constant $0<c_1<\infty$.
On the other hand, conditioning on the time of the last return to 0 before time $l+2$,
and bounding the probability of the latter event by 1, we readily get
$$
\textstyle
\P^{\circ}(X^{\L}_n(l+2)=x, \t^{1}_0< \t^{1}_x)\
\leq (l+2)\P^{\circ}(\t^{0}_x< \t^{0}_0)
=(l+2)\frac{\Q_n(x)}{\Q_n(0)}\P^{\circ}(\t^{x}_0< \t^{x}_x)\,,
\Eq(3.prop4.2)
$$
where the last line follows from reversibility, and where $\Q_n$, defined in
Lemma 2.2 of \cite{BG}, denotes the invariant measure of $X^{\L}_n$. Since
$\P^{\circ}(\t^{x}_0< \t^{x}_x)\leq 1$ we are left to estimate the ratio
of invariant masses in \eqv(3.prop4.2). From the assumption that
$\dist(y,z)=\frac{n}{2}(1-o(1))$, it follows
that $\L_1=n-\L_2=\frac{n}{2}(1-o(1))$. Therefore, by (2.4) of \cite{BG},
$
\frac{\Q_n(x)}{\Q_n(0)}\leq e^{-c_2n}
$
for some constant $0<c_2<\infty$. Gathering our bounds we arrive at
$
p_n^{l+2}(1,z)\leq e^{-c_1n}+(l+2)e^{-c_2n}
$,
which proves the claim of the lemma.
\endproof

\bigskip


\chap{4. Preparations to the verification of Conditions (A1) and (A2).}4

We now capitalize on the estimates of Section 3 and, as a first step towards the verification of Conditions (A1) and (A2),
prove that
these conditions can be replaced by simple ones, where all quantities depending on
the jump chain have been averaged out.
To state our central result we need a little notation.
Set $k_n(t):=\lfloor a_n t\rfloor$. Let $\pi_n^{J,t}(x)$ denote the average number
of visits of $J_n$ to $x$ during the first $k_n(t)$ steps,
$$
\pi_n^{J,t}(x)={k^{-1}_n(t)}\sum_{j=1}^{k_n(t)}\1_{\{J_n(j-1)=x\}}\,,\quad x\in\VV_n\,.
\Eq(4.7)
$$
For $y\in\VV_n$ and $u>0$ further set
$$
h^{u}_n(y)=\sum_{x\in\VV_n}p_n(y,x)\exp\{-uc_n\l_n(x)\}\,,
\Eq(4.8)
$$
and define
$$
\eqalign{
\nu_n^{J,t}(u,\infty)
&= \sum_{j=1}^{k_n(t)}h^{u}_n(J_n(j-1))
=k_n(t)\sum_{y\in\VV_n}\pi_n^{J,t}(y)h^{u}_n(y)\,,
\cr
(\s_n^{J,t})^2(u,\infty)
&=
\sum_{j=1}^{k_n(t)}\left(h^{u}_n(J_n(j-1))\right)^2
=
k_n(t)\sum_{y\in\VV_n}\pi_n^{J,t}(y)\left(h^{u}_n(y)\right)^2\,.
}
\Eq(4.4)
$$
By assumption, the initial distribution, $\mu_n$, is the invariant measure $\pi_n$ of $J_n$.
This implies that the chain variables $(J_n(j), j\geq 1)$ satisfy
$$
P_{\pi_n}(J_n(j)=x)=\pi_n(x)=2^{-n}
\text{for all $x\in\VV_n$, and all $j\geq 1$.}
\Eq(4.6)
$$
Hence
$$
\eqalign{
E_{\pi_n}\left[\pi_n^{J,t}(y)\right]&=\pi_n(y)\,,
\cr
E_{\pi_n}\left[\nu_n^{J,t}(u,\infty)\right]&=\frac{k_n(t)}{a_n}\nu_n(u,\infty)\,,
\cr
E_{\pi_n}\left[(\s_n^{J,t})^2(u,\infty)\right]&=\frac{k_n(t)}{a_n}\s_n^2(u,\infty)\,,
}
\Eq(4.9)
$$
where
$$
\eqalign{
\nu_n(u,\infty)
&=a_n\sum_{x\in\VV_n}\pi_n(x)h^{u}_n(x)\,,
\cr
\s_n^2(u,\infty)
&=a_n\sum_{x\in\VV_n}\pi_n(x)\left(h^{u}_n(x)\right)^2\,.
}
\Eq(4.5)
$$

We will often refer to Proposition \thv(4.prop1) below as to an ergodic theorem.

\proposition{\TH(4.prop1)} {\it
Let $\rho_n>0$ be a decreasing sequence satisfying $\rho_n\downarrow 0$ as $n\uparrow\infty$.
There exists a sequence of subsets $\O^{\tau}_{n,0}\subset\O^{\tau}$ with
%
%
$
\P\left((\O^{\tau}_{n,0})^c\right)<\frac{\theta_nk_n(t)}{\rho_na^2_n}\,,
$
and such that, on $\O^{\tau}_{n,0}$, the following holds: for all $t>0$ and all $u>0$,
$$
\eqalign{
P_{\pi_n}\left(\left|\nu_n^{J,t}(u,\infty)-E_{\pi_n}\left[\nu_n^{J,t}(u,\infty)\right]\right|\geq\e\right)
\leq
\e^{-2}\Theta_n(t,u)\,,\quad\forall\e>0\,,
}
\Eq(4.prop1.1)
$$
where, for some constant $0<c<\infty$,
$$
\Theta_n(t,u)=
\left(\frac{k_n(t)}{a_n}\right)^2\frac{\nu_n^2(u,\infty)}{2^{n}}
+\frac{k_n(t)}{a_n}\s_n^2(u,\infty)
+c\frac{\nu_n(2u,\infty)}{ n^{2}}
+\rho_n\left[\E\nu_n(u,\infty)\right]^2\,.
\Eq(4.prop1.2)
$$
In addition, for all $t>0$ and all $u>0$,
$$
P_{\pi_n}\left((\s_n^{J,t})^2(u,\infty)\geq\e'\right)\leq\frac{k_n(t)}{\e'\, a_n}\s_n^2(u,\infty)\,,\quad\forall\e'>0\,.
\Eq(4.prop1'.2)
$$
}

\remark  Clearly that Proposition \thv(4.prop1) will be useful only if $a_n\gg\theta_n$,
namely, if the number of steps $k_n(t)$ taken by the jump chain is large compared to its stationary time.
From the relations $a_n\sim b_n$ and \eqv(1.1') we see that
if $a_n$ is bounded below then so is $r_n$.
This explains our lower bound on $r_n$ in Assertion (i) of Theorem \thv(1.theo1).

\remark  The small simplification that the choice $\mu_n=\pi_n$ introduces is that
$E_{\mu_n}\left[\pi_n^{J,t}(y)\right]$ exactly is the invariant mass $\pi_n(y)$,
and this in turn yields exact expressions for $E_{\mu_n}\left[\nu_n^{J,t}(u,\infty)\right]$
and $E_{\mu_n}\left[(\s_n^{J,t})^2(u,\infty)\right]$ in \eqv(4.9). When $\mu_n$ is not
the invariant measure of the jump chain one uses first Lemma \thv(3.lemma1) to approximate
$E_{\mu_n}\left[\pi_n^{J,t}(y)\right]$ by its invariant mass.
Proposition \thv(4.prop1) is then proved along the same lines.
Choosing $\pi_n$  for initial distribution is generic inasmuch
as the error introduced by this approximation is typically  negligible.

\remark  A slightly different version of Proposition \thv(4.prop1) will be needed to deal
with the case of extreme scales, where a representation of the landscape will be substituted for
the original one. See Proposition \thv(6.prop1) of Section 6.2.

\proofof{Proposition \thv(4.prop1)} The upper bound \eqv(4.prop1'.2) plainly results from a first order
Tchebychev inequality and the expression \eqv(4.9) of $E_{\pi_n}\left[(\s_n^{J,t})^2(u,\infty)\right]$.
The proof of \eqv(4.prop1.1) is a little more involved.
Using a second order Tchebychev inequality together with the expressions
\eqv(4.9) and \eqv(4.5) of $E_{\pi_n}\left[\nu_n^{J,t}(u,\infty)\right]$ and $\nu_n(u,\infty)$,  we get,
$$
P_{\pi_n}\left(\left|\nu_n^{J,t}(u,\infty)-E_{\pi_n}\left[\nu_n^{J,t}(u,\infty)\right]\right|\geq\e\right)
\leq\e^{-2}
E_{\pi_n}\Bigl[
k_n(t)\sum_{y\in\VV_n}\left(\pi_n^{J,t}(y)-\pi_n(y)\right)h^{u}_n(y)
\Bigr]^2\,.
\Eq(4.prop1.6)
$$
Expanding the r.h.s. of \eqv(4.prop1.6) yields
$$
\e^{-2}
k^2_n(t)\sum_{x\in\VV_n}\sum_{y\in\VV_n}h^{u}_n(x)h^{u}_n(y)E_{\pi_n}\left(\pi_n^{J,t}(x)-\pi_n(x)\right)\left(\pi_n^{J,t}(y)-\pi_n(y)\right)\,.
\Eq(4.prop1.7)
$$
In view of \eqv(4.7), setting
$
\D_{ij}(x,y)=P_{\pi_n}\left(J_n(i-1)=x, J_n(j-1)=y\right)-\pi_n(x)\pi_n(y)
$,
the expectation in \eqv(4.prop1.7) may be expressed as
$$
k^2_n(t)E_{\pi_n}\left(\pi_n^{J,t}(x)-\pi_n(x)\right)\left(\pi_n^{J,t}(y)-\pi_n(y)\right)
=\sum_{i=1}^{k_n(t)}\sum_{j=1}^{k_n(t)}\D_{ij}(x,y)\,.
\Eq(4.prop1.9)
$$
For $\theta_n$ defined as in Proposition \thv(3.prop1), we now break the sum in the r.h.s. of \eqv(4.prop1.9) into three terms:
$$
\eqalign{
(\overline I)&=2\sum_{1\leq i\leq k_n(t)}\sum_{i+\theta_n\leq j\leq k_n(t)}\D_{ij}(x,y)\,,
\cr
(\overline{II})&=\sum_{1\leq i\leq k_n(t)}\1_{\{i=j\}}\D_{ij}(x,y)\,,
\cr
(\overline{III})&=2\sum_{1\leq i\leq k_n(t)}\sum_{i<j<i+\theta_n}\D_{ij}(x,y)\,.
}
\Eq(4.prop1.10)
$$
Consider first $(\overline I)$. By Proposition \thv(3.prop1),
$$
\eqalign{
(\overline I)
&\leq
2\sum_{1\leq i\leq k_n(t)}\sum_{\lfloor(i+\theta_n)/2\rfloor\leq l\leq \lceil k_n(t)/2\rceil}[\D_{i(2l)}(x,y)+\D_{i(2l+1)}(x,y)]
\cr
&\leq
2\sum_{1\leq i\leq k_n(t)}\sum_{\lfloor(i+\theta_n)/2\rfloor\leq l\leq \lceil k_n(t)/2\rceil}2|\d_n|\pi_n(x)\pi_n(y)
\cr
&\leq
|\d_n|k^2_n(t)\pi_n(x)\pi_n(y)\,.
\cr
}
\Eq(4.prop1.10')
$$
where $|\d_n|\leq 2^{-n}$.
Turning to the term $(\overline{II})$, we have,
$$
\eqalign{
(\overline{II})
&=\sum_{1\leq i\leq k_n(t)}\D_{ii}(x,x)\1_{\{x=y\}}
\cr
&=k_n(t)\left[P_{\pi_n}\left(J_n(i-1)=x\right)-\pi^2_n(x)\right]\1_{\{x=y\}}
\cr
&=k_n(t)\pi_n(x)(1-\pi_n(x))\1_{\{x=y\}}\,,
\cr
}
\Eq(4.prop1.11)
$$
where the last equality follows from \eqv(4.6). Finally,
$$
\eqalign{
(\overline{III})
&\leq
2\sum_{i=1}^{k_n(t)}\sum_{l=1}^{\theta_n-1}P_{\pi_n}\left(J_n(i-1)=x, J_n(i+l-1)=y\right)
\cr
&\leq 2\sum_{i=1}^{k_n(t)}\sum_{l=1}^{\theta_n-1}
P_{\pi_n}\left(J_n(i-1)=x\right)P_{\pi_n}\left(J_n(i+l-1)=y\mid J_n(i-1)=x\right)
\cr
&= 2k_n(t)\pi_n(x)\sum_{l=1}^{\theta_n-1}p_n^{l}(x,y)\,,
}
\Eq(4.prop1.12)
$$
where $p_n^{l}(\cdot,\cdot)$ denote the $l$-steps transition matrix of $J_n$.
Combining our bounds on $(\overline{I}),(\overline{II})$, and $(\overline{III})$ with \eqv(4.prop1.7) we get that, for all $\e>0$,
$$
P_{\pi_n}\left(\left|\nu_n^{J,t}(u,\infty)-E_{\pi_n}\left[\nu_n^{J,t}(u,\infty)\right]\right|\geq\e\right)
\leq\e^{-2}[(I)+(II)+(III)]\,,
\Eq(4.prop1.13)
$$
where
$$
\eqalign{
(I)
&=
|\d_n|k^2_n(t)\sum_{x\in\VV_n}\sum_{y\in\VV_n}h^{u}_n(x)h^{u}_n(y)(t)\pi_n(x)\pi_n(y)\,,
\cr
(II)&=
k_n(t)\sum_{x\in\VV_n}\sum_{y\in\VV_n}h^{u}_n(x)h^{u}_n(y)\pi_n(x)(1-\pi_n(x))\1_{\{x=y\}}\,,
\cr
(III)
&=2k_n(t)\sum_{x\in\VV_n}\sum_{y\in\VV_n}h^{u}_n(x)h^{u}_n(y)\pi_n(x)\sum_{l=1}^{\theta_n-1}p_n^{l}(x,y)\,.
}
\Eq(4.prop1.14)
$$
By \eqv(4.5),
$$
\eqalign{
(I)&\leq
\left(\frac{k_n(t)}{a_n}\right)^2\frac{\nu_n^2(u,\infty)}{2^{n}}\,,
\cr
(II)&\leq
\frac{k_n(t)}{a_n}\s_n^2(u,\infty)\,.
}
\Eq(4.prop1.15)
$$
To further express the third term in \eqv(4.prop1.14) note that, by \eqv(4.8),
$$
\sum_{y\in\VV_n}p_n^{l}(x,y)h^{u}_n(y)
=\sum_{y\in\VV_n}p_n^{l}(x,y)\sum_{z\in\VV_n}p_n(y,z)e^{-uc_n\l_n(z)}
=\sum_{z\in\VV_n}p_n^{l+1}(x,z)e^{-uc_n\l_n(z)}\,,
\Eq(4.prop1.16)
$$
and,
$$
\eqalign{
\sum_{x\in\VV_n}\pi_n(x)h^{u}_n(x)p_n^{l+1}(x,z)
=&
\sum_{y\in\VV_n}e^{-uc_n\l_n(y)}\sum_{x\in\VV_n}\pi_n(x)p_n(x,y)p_n^{l+1}(x,z)
\cr
=&
\sum_{y\in\VV_n}e^{-uc_n\l_n(y)}\pi_n(y)p_n^{l+2}(y,z)\,,
}
\Eq(4.prop1.17)
$$
where the last equality follows by reversibility.
Hence,
$$
\eqalign{
(III)
=&2k_n(t)\sum_{l=1}^{\theta_n-1}\sum_{z\in\VV_n}\left[\sum_{x\in\VV_n}\pi_n(x)h^{u}_n(x)p_n^{l+1}(x,z)\right]e^{-uc_n\l_n(z)}\,,
\cr
=&2\sum_{l=1}^{\theta_n-1}k_n(t)\sum_{z\in\VV_n}\sum_{y\in\VV_n}\pi_n(y)e^{-uc_n(\l_n(y)+\l_n(z))}p_n^{l+2}(y,z)
\cr
=&2\sum_{l=1}^{\theta_n-1}[(III)_{1,l}+(III)_{2,l}]
\,.
}\Eq(4.prop1.18)
$$
where, distinguishing the cases $z=y$ and $z\neq y$,
$$
\eqalign{
(III)_{1,l}
&=\sum_{z\in\VV_n}k_n(t)\pi_n(z)e^{-2uc_n\l_n(z)}p_n^{l+2}(z,z)\,,
\cr
(III)_{2,l}
&=\sum_{z\in\VV_n}\sum_{y\in\VV_n : y\neq z}k_n(t)\pi_n(y)e^{-uc_n(\l_n(y)+\l_n(z))}p_n^{l+2}(y,z)\,.
}
\Eq(4.prop1.19)
$$
One easily checks that $\theta_n\leq 2m$ with $m\leq n^2$. Thus, by Proposition \thv(3.prop3),
$$
\sum_{l=1}^{\theta_n-1}(III)_{1,l}
=\sum_{z\in\VV_n}k_n(t)\pi_n(z)e^{-2uc_n\l_n(z)}\sum_{l=1}^{\theta_n-1}p_n^{l+2}(z,z)
\leq c n^{-2}\nu_n(2u,\infty)\,.
\Eq(4.prop1.20)
$$
for some constant $0<c<\infty$.

The next lemma is designed to deal with the second sum in the last line of \eqv(4.prop1.18).

\lemma{\TH(4.lemma1)} {\it
Let $\rho_n>0$ be a decreasing sequence satisfying $\rho_n\downarrow 0$ as $n\uparrow\infty$.
There exists a sequence of subsets $\O^{\tau}_{n,0}\subset\O^{\tau}$ with
$
\P\left((\O^{\tau}_{n,0})^c\right)<\frac{\theta_nk_n(t)}{\rho_na^2_n}\,,
$
and such that, on $\O^{\tau}_{n,0}$,
$$
\sum_{l=1}^{\theta_n-1}(III)_{2,l}<\rho_n\left[\E\nu_n(u,\infty)\right]^2\,.
\Eq(4.lem1.2)
$$
}

\proof  By a first order Tchebychev inequality, for all $\eta>0$,
$
\P\left(\sum_{l=1}^{\theta_n-1}(III)_{2,l}\geq\eta\right)
\leq
\eta^{-1}\sum_{l=1}^{\theta_n-1}\E(III)_{2,l}
$.
Next,  for all $y\neq z\in\VV_n\times \VV_n$, by independence,
$
\E\left[e^{-uc_n(\l_n(y)+\l_n(z))}\right]=\left[a^{-1}_n\E\nu_n(u,\infty)\right]^2
$.
Thus,
$$
\sum_{l=1}^{\theta_n-1}\E(III)_{2,l}
\leq
\frac{k_n(t)}{a^2_n}\left[\E\nu_n(u,\infty)\right]^2\sum_{l=1}^{\theta_n-1}\sum_{z\in\VV_n}p_n^{l+2}(y,z)
\leq
\frac{\theta_nk_n(t)}{a^2_n}\left[\E\nu_n(u,\infty)\right]^2\,,
$$
yielding
$
\P\left(\sum_{l=1}^{\theta_n-1}(III)_{2,l}\geq\eta\right)
\leq\frac{\theta_nk_n(t)}{\eta a^2_n}\left[\E\nu_n(u,\infty)\right]^2
$.
The lemma now easily follows. \endproof

Collecting the bounds \eqv(4.prop1.15), \eqv(4.prop1.20), and \eqv(4.lem1.2),
and combining them with \eqv(4.prop1.13), we obtain that under the assumptions
and with the notations of Lemma \thv(4.lemma1), on $\O^{\tau}_{n,0}$,
for all $t>0$ and all $u>0$,
$$
\eqalign{
&P_{\pi_n}\left(\left|\nu_n^{J,t}(u,\infty)-E_{\pi_n}\left[\nu_n^{J,t}(u,\infty)\right]\right|\geq\e\right)
\cr
\leq&\e^{-2}\left\{
\left(\frac{k_n(t)}{a_n}\right)^2\frac{\nu_n^2(u,\infty)}{2^{n}}
+\frac{k_n(t)}{a_n}\s_n^2(u,\infty)
+c\frac{\nu_n(2u,\infty)}{ n^{2}}
+\rho_n\left[\E_{\pi_n}\nu_n(u,\infty)\right]^2
\right\}\,.
}
\Eq(4.prop1.21)
$$
for some constant $0<c<\infty$.
The proof of Proposition \thv(4.prop1) is done. \endproof

\bigskip


\chap{5. Intermediate and short scales.}5

Set
$$
\g_n(x)=r_n^{-1}\t_n(x)\,,\quad x\in\VV_n\,.
\Eq(5.1)
$$
We call $(\g_n(x), x\in\VV_n)$  the re-scaled landscape.
With this notation, and in the present setting, the quantities $\nu_n$ and $\s_n^2$ of \eqv(4.5) reads
$$
\nu_n(u,\infty)=
\frac{a_n}{2^n}\sum_{x\in\VV_n}
e^{-u/\g_n(x)}\,,
\Eq(5.2)
$$
$$
\s_n^2(u,\infty)
=
\frac{a_n}{2^n}\sum_{x,\in\VV_n}\sum_{x'\in\VV_n}
p^{2}_n(x,x')e^{-u/\g_n(x)}e^{-u/\g_n(x')}\,,
\Eq(5.3)
$$
where $p^{2}_n(\cdot,\cdot)$ are the two steps transition probabilities of $J_n$.

\bigskip
\line{\bf 5.1 Chain independent estimates.\hfill}

In the following two propositions (Proposition \thv(5.prop6) and Proposition \thv(5.prop7))
we collect the chain independent results
needed to establish the validity of Conditions (A1), (A2), and (A3).
Note that Condition (A0') reads
$$
{\nu_n(u,\infty)}/{a_n}=o(1)\,.
\Eq(1.A0'bis)
$$
Now this will hold true as a by-product of our convergence results for $\nu_n(u,\infty)$.

\proposition{\TH(5.prop6)} [Intermediate space scales] {\it  Let $r_n$ be an intermediate space scale and
choose $a_n\sim b_n$ in \eqv(5.2) and \eqv(5.3).
Let $\nu^{int}$ be defined in \eqv(1.prop1.1) and assume that  $\b\geq\b_c(\varepsilon)$.

\item{i)} If $\sum_{n}a_n/2^n<\infty$ then there exists a subset $\O^{\tau}_1\subset\O^{\tau}$
with $\P(\O^{\tau}_1)=1$ such that, on $\O^{\tau}_1$, the following holds: for all $u>0$
$$
\eqalign{
&\lim_{n\rightarrow\infty}\nu_n(u,\infty)=\nu^{int}(u,\infty)\,,
\cr
&\lim_{n\rightarrow\infty}n\s_n^2(u,\infty)=\nu^{int}(2u,\infty)\,,
}
\Eq(5.prop6.1)
$$
and
$
\lim_{\d\rightarrow 0}\lim_{n\rightarrow\infty}\int_0^{\d}\nu_n(u,\infty)=0
$.

\item{ii)} If $\sum_{n}a_n/2^n=\infty$ then there exists a
sequence of subsets $\O^{\tau}_{1,n}\subset\O^{\tau}$
with $\P(\O^{\tau}_{1,n})\geq 1-o(1)$ and such that for all $n$ large enough, on $\O^{\tau}_{1,n}$,
the following holds: for all $u>0$
$$
\eqalign{
&
\left|\nu_{n}(u,\infty)-\E[\nu_{n}(u,\infty)]\right|
\leq  2(a_n/2^n)^{1/4}\bigl(1 \vee\sqrt{2\nu^{int}(2u,\infty)}\bigr)\,,
\cr
&
\left|n\s^2_{n}(u,\infty)-n\E[\s^2_{n}(u,\infty)]\right|
\leq (a_n/2^n)^{1/4}\bigl(1 \vee\sqrt{2\nu^{int}(u,\infty)}\bigr)\,,
\cr
}
\Eq(5.prop6.2)
$$
where
$$
\eqalign{
&
\lim_{n\rightarrow\infty}\E[\nu_{n}(u,\infty)]=\nu^{int}(u,\infty)\,,
\cr
&
\lim_{n\rightarrow\infty}n\E[\s^2_{n}(u,\infty)]=\nu^{int}(2u,\infty)\,,
\cr
}
\Eq(5.prop6.2bis)
$$
and, for all $0<\d\leq \d_0$ and small enough $\d_0$,
$
\int_0^{\d}\nu_n(u,\infty)du\leq 
c_0\d^{1-\a(\varepsilon)}\sfrac{\a(\varepsilon)}{1-\a(\varepsilon)}\G(\a(\varepsilon))\,,
$
for some numerical constant $0<c_0<\infty$.
}


\proposition{\TH(5.prop7)} [Short space scales] {\it Let $r_n$ be an intermediate space scale
and choose $a_n\sim b_n$ in \eqv(5.2) and \eqv(5.3). There exists a subset
$\O^{\tau}_1\subset\O^{\tau}$ with $\P(\O^{\tau}_1)=1$ such that, on $\O^{\tau}_1$,
the following holds:
$
\lim_{\d\rightarrow 0}\lim_{n\rightarrow\infty}\int_0^{\d}\nu_n(u,\infty)=0
$
and:
\item{i)} if $n/a_n=o(1)$ then for all $u>0$,
$$
\eqalign{
&\lim_{n\rightarrow\infty}\nu_n(u,\infty)=1\,,
\cr
&\lim_{n\rightarrow\infty}n\s_n^2(u,\infty)=1\,;
}
\Eq(5.prop7.1)
$$
\item{ii)} if $a_n/n=\OO(1)$  then for all $u>0$, setting  $C=1+\lim_{n\rightarrow\infty}a_n/n$,
$$
\eqalign{
&\lim_{n\rightarrow\infty}\nu_n(u,\infty)=1\,,
\cr
&\lim_{n\rightarrow\infty}a_n\s_n^2(u,\infty)=C\,.
}
\Eq(5.prop7.2)
$$
}

The proofs of these propositions, which are given at the end of this subsection, rely on the following three key lemmata.

\lemma{\TH(5.lemma7)} {\it Let $r_n$ be an intermediate or short space scale. Then
for all $0\leq\varepsilon\leq 1$ and all  $0<\beta<\infty$ satisfying $\b\geq \b_c(\varepsilon)$,
choosing $a_n\sim b_n$ in \eqv(5.2),
$$
\lim_{n\rightarrow\infty}\E[\nu_n(u,\infty)]
=
\cases
\nu^{int}(u,\infty)
&\hbox{if $\quad 0<\a(\varepsilon)\leq 1$} , \,\,\,\cr
1
, &\hbox{if $\quad\a(\varepsilon)=0$.}\,\,\,\cr
\endcases
\,,\quad u>0\,.
\Eq(5.lem7.1)
$$
where $\nu^{int}$ is defined in \eqv(1.prop1.1). Moreover, for all $L\geq 0$ such that $a_nL/2^n=o(1)$,
$$
\P\left(\left|\nu_n(u,\infty)-\E[\nu_n(u,\infty)]\right|
\geq 2\sqrt{a_nL/2^n}\sqrt{\E[\nu_n(2u,\infty)]}\right)
\leq e^{-L}\,.
\Eq(5.lem7.2)
$$
}

\lemma{\TH(5.lemma8')} {\it
$
\displaystyle
\E[\s_n^2(u,\infty)]
=\frac{\E[\nu_n(2u,\infty)]}{n}+\frac{(\E[\nu_n(u,\infty)])^2}{a_n}\frac{n-1}{n}\,.
$}

\lemma{\TH(5.lemma8)} {\it
Let $r_n$ be an intermediate or short space scale.
Then for all $0\leq\varepsilon\leq 1$ and all  $0<\beta<\infty$ satisfying $\b\geq \b_c(\varepsilon)$,
choosing $a_n\sim b_n$ in \eqv(5.3), the
following holds for all $L>0$:

\item{(i)} if ${n}/{a_n}=o(1)$,
$$
\P\left(\left|\s_n^2(u,\infty)-\E[\s_n^2(u,\infty)]\right|
\geq n^{-1}\sqrt{a_nL/2^n}\sqrt{\E\left[\nu_n(u,\infty)\right]}\right)
\leq L^{-1}\,,
\Eq(5.lem8.2)
$$

\item{(ii)} if $a_n/n=\OO(1)$,
$$
\P\left(\left|\s_n^2(u,\infty)-\E[\s_n^2(u,\infty)]\right|
\geq \frac{n}{a_n}\sqrt{a_nL/2^n}\sqrt{\overline C\E\left[\nu_n(u,\infty)\right]}\right)
\leq L^{-1}\,,
\Eq(5.lem8.2bis)
$$
for some constant $0<\overline C<\infty$.
}

\proofof{Lemma \thv(5.lemma7)} We first prove \eqv(5.lem7.1).
By \eqv(5.2), $\E[\nu_n(u,\infty)]=a_n\E e^{-u/\g_n(0)}$, $0\in\VV_n$.
For fixed $u>0$ set $f(y)=e^{-u/y}$.
Thus $f(0)=0$, $f'(y)=(u/y^2)e^{-u/y}$ and, integrating by part,
$$
\E[\nu_n(u,\infty)]
=
a_n\int_{0}^{\infty}f'(y)\P\left(\g_n(0)\geq y\right)dy
=
\frac{a_n}{b_n}\int_{0}^{\infty}f'(y)h_n(y)dy\,.
\Eq(5.lem7.3)
$$
Set  $I_n(a,b)=\int_{a}^{b}f'(y)h_n(y)dy$, $a\leq b$. Given  $0<\hat\zeta<1$ and $\zeta>1$,
we may rewrite \eqv(5.lem7.3) as
$$
\E[\nu_n(u,\infty)]
=
(1+o(1))\left[I_n\bigl(0,r_n^{-1/2}\bigr)+I_n\bigl(r_n^{-1/2},\hat\zeta\bigr)+I_n\bigl(\hat\zeta,\zeta\bigr)+I_n(\zeta,\infty)\right]\,,
\Eq(5.lem7.3')
$$
where we used the  assumption that $a_n/b_n\sim 1$.
We will now show that,  as $n\rightarrow\infty$, for small enough $\hat\zeta$ and large enough $\zeta$,
the leading contribution to \eqv(5.lem7.3') comes from $I_n\bigl(\hat\zeta,\zeta\bigr)$.
To do so we first use that by \eqv(1.1') and the rough upper bound
$
h_n(y)\leq b_n
$,
$
\textstyle
I_n\bigl(0,r_n^{-1/2}\bigr)\leq b_n\int_{0}^{1/\sqrt{r_n}}f'(y)dy={e^{-u\sqrt{r_n}}}/{\P(\t_n(x)\geq r_n)}
$,
and, together with the gaussian tail estimates \eqv(2.11'),  this entails
$$
\lim_{n\rightarrow\infty}I_n\bigl(0,r_n^{-1/2}\bigr)=0.
\Eq(5.lem7.3'')
$$
Next, by Lemma \thv(2.lemma6), (ii), with $\d=1/2$,
$
I_n\bigl(r_n^{-1/2},\hat\zeta\bigr)
\leq 2(1+o(1))\int_{0}^{\hat\zeta}f'(y)y^{-(3/4)\a_n}dy
$
for all $0<\hat\zeta<1$, where $0\leq \a_n=\a(\varepsilon)+o(1)$.
Now, there exists $\zeta^*\equiv\zeta^*(u)>0$ such that, for all $\hat\zeta<\zeta^*$,
$f'(y)y^{-(3/4)\a_n}$ is strictly increasing on $[0,\hat\zeta]$. Hence,
for all $\hat\zeta<\min(1,\zeta^*)$,
$
I_n\bigl(r_n^{-1/2},\hat\zeta\bigr)
\leq 2(1+o(1))u\hat\zeta^{-1+(3/4)[\a(\varepsilon)+o(1)]}e^{-{u}/{\hat\zeta}}
$,
implying that
$$
\lim_{n\rightarrow\infty}I_n\bigl(r_n^{-1/2},\hat\zeta\bigr)
\leq 2u\hat\zeta^{-1+(3/4)\a(\varepsilon)}e^{-{u}/{\hat\zeta}}\,,\quad \hat\zeta<\min(1,\zeta^*)\,.
\Eq(5.lem7.3''')
$$
To deal with $I_n\bigl(\hat\zeta,\zeta\bigr)$ note that
by Lemma \thv(2.lemma6), (i), $h_n(y)\rightarrow y^{-\a(\varepsilon)}$, $n\rightarrow\infty$,
where the convergence is uniform in  $\hat\zeta\leq y\leq \zeta$ since, for each $n$,
$h_n(y)$ is a monotone function, and since the limit, $y^{-\a(\varepsilon)}$,
is continuous. Hence,
$$
\lim_{n\rightarrow\infty}I_n\bigl(\hat\zeta,\zeta\bigr)
=\lim_{n\rightarrow\infty}\int_{\hat\zeta}^{\zeta}f'(y)h_n(y)dy
=\int_{\hat\zeta}^{\zeta}f'(y)y^{-\a(\varepsilon)}dy\,.
\Eq(5.lem7.4)
$$
It remains to bound $I_n(\zeta,\infty)$. By
\eqv(2.lem6.1) of Lemma \thv(2.lemma6),
$
I_n(\zeta,\infty)
=\int_{\zeta}^{\infty}f'(y)h_n(y)dy
=(1+o(1))\int_{\zeta}^{\infty}f'(y)y^{-\a_n}dy\,,
$
where again $0\leq \a_n=\a(\varepsilon)+o(1)$.
Thus, for $0<\d<1$ arbitrary we have, taking $n$ large enough, that for all $y\geq\zeta>1$,
$
f'(y)y^{-\a_n}
\leq f'(y)y^{-\a(\varepsilon)+\d}
\leq u/y^{2-\d}
$.
Therefore
$
I_n(\zeta,\infty)
\leq(1+o(1))\frac{1}{1-\d}\zeta^{-(1-\d)}
$,
and, choosing e.g\. $\d=1/2$,
$$
\lim_{n\rightarrow\infty}I_n(\zeta,\infty)\leq 2u\zeta^{-1/2}\,.
\Eq(5.lem7.5)
$$
Collecting \eqv(5.lem7.3''), \eqv(5.lem7.3'''), \eqv(5.lem7.4) and \eqv(5.lem7.5),
we obtain that for all $\zeta>1$ and $\hat\zeta<\min(1,\zeta^*)$,
$$
\lim_{n\rightarrow\infty}\E[\nu_n(u,\infty)]
=\int_{\hat\zeta}^{\zeta}f'(y)y^{-\a(\varepsilon)}dy
+\RR(\hat\zeta,\zeta)\,,
\Eq(5.lem7.6)
$$
where
$
0\leq \RR(\hat\zeta,\zeta)\leq
+2u\hat\zeta^{-1+(3/4)\a(\varepsilon)}e^{-{u}/{\hat\zeta}}
+2u\zeta^{-1/2}
$.
Finally, passing to the limit $\hat\zeta\rightarrow 0$ and $\zeta\rightarrow\infty$ in \eqv(5.lem7.6) yields
$$
\lim_{n\rightarrow\infty}\E[\nu_n(u,\infty)]=\int_{0}^{\infty}f'(y)y^{-\a(\varepsilon)}dy\,,
\Eq(5.lem7.8)
$$
which is (in particular)
%
%
valid for all $0\leq\a(\varepsilon)\leq 1$. For $\a(\varepsilon)=0$,
$
\int_{0}^{\infty}f'(y)y^{-\a(\varepsilon)}dy=\int_{0}^{\infty}f'(y)dy=1
$
while for $0<\a(\varepsilon)\leq 1$,
$
\int_{0}^{\infty}f'(y)y^{-\a(\varepsilon)}dy=u^{-\a(\varepsilon)}\a(\varepsilon)\G(\a(\varepsilon))
$.
Thus \eqv(5.lem7.1) is proven.

It remains to prove \eqv(5.lem7.2).
For this we will use Bennett's bound \cite{Ben} for the tail behavior of sums of random variables,
which states that if $(X(x),\, x\in\VV_n)$ is a family of i\.i\.d\. centered random variables
that satisfies $\max_x|X(x)|\leq \bar a$ then, setting
$\tilde b^2=\sum_{x\in\VV_n}\E X^2(x)$,  for all $\bar b^2\geq \tilde b^2$,
$$
\P\Bigl(\Bigl|\sum_{x\in\VV_n}X(x)\Bigr|>t\Bigr)
\leq\exp\left\{
\frac{t}{\bar a}-\left(\frac{t}{\bar a}+\frac{\bar b^2}{\bar a^2}\right)\log\left(1+\frac{\bar at}{\bar b^2}\right)
\right\}\,,\quad t\geq 0\,.
\Eq(5.lem7.10)
$$
The behavior of the r\.h\.s\. of \eqv(5.lem7.10) varies depending on the relative
size of $t$ and of the ratio $\bar b^2/\bar a$. Note in particular that
for $t<\bar b^2/(2\bar a)$, \eqv(5.lem7.10) simplifies to
$$
\textstyle{\P\left(\left|\sum_{x\in\VV_n}X(x)\right|\geq t\right)}
\leq\exp\left\{-{t^2}/{4\bar b^2}\right\}\,.
\Eq(5.lem7.11)
$$
To make use of Bennett's bound we set $X(x)=e^{-u/\g_n(x)}-\E e^{-u/\g_n(x)}$ so that
$$
\P\left(\left|\nu_n(u,\infty)-\E[\nu_n(u,\infty)]\right|\geq \theta\right)
=\textstyle{\P\left(\left|\sum_{x\in\VV_n}X(x)\right|\geq 2^n a_n^{-1}\theta\right)}\,.
\Eq(5.lem7.9)
$$
Since $|X(x)|\leq 1$ and
$
\sum_{x\in\VV_n}\E X^2(x)\leq 2^n\E e^{-2u/\g_n(x)}=2^na_n^{-1}\E[\nu_n(2u,\infty)]
$,
we may choose $\bar a=1$ and $\bar b^2=2^na_n^{-1}\E[\nu_n(2u,\infty)]$.
Then, by  \eqv(5.lem7.11)  and  \eqv(5.lem7.9), for all $L>0$,
$$
\P\left(\left|\nu_n(u,\infty)-\E[\nu_n(u,\infty)]\right|
\geq 2\sqrt{a_nL/2^n}\sqrt{\E[\nu_n(2u,\infty)]}\right)
\leq e^{-L}\,,
\Eq(5.lem7.15)
$$
where we chose $\theta^2= a_n L2^{-n+2}\E[\nu_n(2u,\infty)]$.
This choice is permissible provided that
$\theta\leq \E[\nu_n(2u,\infty)]/2$. In view of \eqv(5.lem7.1)
this will be verified for all $n$ large enough whenever $\theta\downarrow 0$
as $n\uparrow\infty$, i\.e\. whenever $a_nL/2^n=o(1)$. Thus
\eqv(5.lem7.2) is established, and the lemma proven.\endproof

We skip the elementary proof of Lemma \thv(5.lemma8').


\proofof{Lemma \thv(5.lemma8)}
For $u>0$ and $l\geq 1$ set
$$
\s_n^l(u,\infty)=a_n\sum_{y\in\VV_n}\pi_n(y)\left(h^{u}_n(y)\right)^l\,.
\Eq(5.lem8.1)
$$
If for $l=1$, $\s_n^l(u,\infty)=\nu_n(u,\infty)$
is a sum of independent random variables, this is no longer true when $l=2$. In this case
we simply resort to a second order Tchebychev inequality to write
$$
\eqalign{
\P\left(\left|\s_n^l(u,\infty)-\E[\s_n^l(u,\infty)]\right|\geq t\right)
&\leq t^{-2} \E\left[\s_n^l(u,\infty)-\E[\s_n^l(u,\infty)]\right]^2
= t^{-2}[\theta_1+\theta_2],
\cr
}
\Eq(5.lem8.4)
$$
where
$$
\eqalign{
&\theta_1=\left(\frac{a_n}{2^n}\right)^2\sum_{y\in\VV_n}\E\left[(h^{u}_n(y))^l-\E(h^{u}_n(y))^l\right]^2\,,
\cr
&\theta_2=\left(\frac{a_n}{2^n}\right)^2\sum_{{y,y'\in\VV_n\times\VV_n}\atop {y\neq y'}}
\E\left\{\left[(h^{u}_n(y))^l-\E(h^{u}_n(y))^l\right]\left[(h^{u}_n(y'))^l-\E(h^{u}_n(y'))^l\right]\right\}\,.
}
\Eq(5.lem8.5)
$$
On the one hand, we clearly have,
$$
\theta_1=\frac{a_n}{2^n}\E[\s_n^{2l}(u,\infty)]-\frac{1}{2^n}(\E[\s_n^{l}(u,\infty)])^2\leq \frac{a_n}{2^n}\E[\s_n^{2l}(u,\infty)].
\Eq(5.lem8.5')
$$
On the other hand, after some lengthy but simple calculations, we obtain that
$$
\eqalign{
\theta_2
&\leq
\frac{n(n-1)}{2^{n+1}}\Biggl[
\frac{a_n}{n^{2l}}\E\left[\nu_n(u,\infty)\right]
+2\frac{(\E[\nu_n(u,\infty)])^2}{n^l}\left(\frac{\E[\nu_n(u,\infty)]}{a_n}+\frac{2}{n}\right)^{l-1}
\cr
&+\frac{1}{a_n}(\E[\nu_n(u,\infty)])^3\left(\frac{\E[\nu_n(u,\infty)]}{a_n}+\frac{1}{n}\right)^{2(l-1)}
\Biggr]
\,.
}
\Eq(5.lem8.20)
$$

We now specialize the bounds \eqv(5.lem8.5') and \eqv(5.lem8.20)
on $\theta_1$ and $\theta_2$ according to whether
${n}/{a_n}=o(1)$ or ${a_n}/n=\OO(1)$.
The resulting bounds will be valid under the
assumptions of Lemma \thv(5.lemma7), and for large enough  $n$.
Assume first that ${n}/{a_n}=o(1)$.
Then
$$
\eqalign{
\theta_2
&\leq
\frac{\E\left[\nu_n(u,\infty)\right]}{2n^{2(l-1)}}\frac{a_n}{2^{n}}\Biggl[1
+2(\E[\nu_n(u,\infty)])^2\frac{n}{a_n}\left(2+\frac{n}{a_n}\E[\nu_n(u,\infty)]\right)^{l-1}
\cr
&+\frac{1}{a_n}(\E[\nu_n(u,\infty)])^3\left(\frac{n}{a_n}\right)^2\left(2+\frac{n}{a_n}\E[\nu_n(u,\infty)]\right)^{2(l-1)}
\Biggr]
\,.
}
\Eq(5.lem8.21)
$$
Thus, in view of \eqv(5.lem7.1), for all $n$ large enough,
$$
\theta_2\leq \frac{(1+o(1))}{2}\frac{\E\left[\nu_n(u,\infty)\right]}{n^{2(l-1)}}\frac{a_n}{2^{n}}\,.
\Eq(5.lem8.22)
$$
We prove in a similar way that, for all $n$ large enough,
$$
\theta_1\leq \frac{(1+o(1))}{n}\frac{\E\left[\nu_n(u,\infty)\right]}{n^{2(l-1)}}\frac{a_n}{2^{n}}\,.
\Eq(5.lem8.23)
$$
From the last two bounds it follows that, for all $n$ large enough,
$$
\theta_1+\theta_2\leq\frac{\E\left[\nu_n(u,\infty)\right]}{n^{2(l-1)}}\frac{a_n}{2^{n}}\,.
\Eq(5.lem8.24)
$$

Assume now that $a_n/n=\OO(1)$ and $a_n\geq 1$. Reasoning as above it easily follows from \eqv(5.lem8.20)
and  \eqv(5.lem8.5') respectively that there exist constants $0<\overline C, \overline C'<\infty$
such that, for all  $n$ large enough,
$$
\theta_2\leq \overline C\E\left[\nu_n(u,\infty)\right]\frac{a_n}{2^{n+1}}\left(\frac{n}{a_n}\right)^2\,,
\Eq(5.lem8.25)
$$
and
$$
\theta_1\leq \overline C'\frac{a_n}{2^n}E\left[\nu_n(u,\infty)\right]\,,
\Eq(5.lem8.26)
$$
so that
$$
\theta_1+\theta_2\leq \overline C\E\left[\nu_n(u,\infty)\right]\frac{a_n}{2^{n}}\left(\frac{n}{a_n}\right)^2\,.
\Eq(5.lem8.27)
$$

Inserting \eqv(5.lem8.24) in \eqv(5.lem8.4) and choosing
$
t=n^{-(l-1)}\sqrt{\frac{a_nL}{2^n}\E\left[\nu_n(u,\infty)\right]}
$
yields
$$
\P\left(\left|\s_n^l(u,\infty)-\E[\s_n^l(u,\infty)]\right|
\geq n^{-(l-1)}\sqrt{a_nL/2^n}\sqrt{\E\left[\nu_n(u,\infty)\right]}\right)
\leq L^{-1}\,.
\Eq(5.lem8.2')
$$
Similarly it follows  from \eqv(5.lem8.27)
and the choice
$
t=\frac{n}{a_n}\sqrt{\frac{a_nL}{2^n}\overline C\E\left[\nu_n(u,\infty)\right]}
$
that
$$
\P\left(\left|\s_n^l(u,\infty)-\E[\s_n^l(u,\infty)]\right|
\geq \frac{n}{a_n}\sqrt{a_nL/2^n}\sqrt{\overline C\E\left[\nu_n(u,\infty)\right]}\right)
\leq L^{-1}\,.
\Eq(5.lem8.2'bis)
$$
Taking $l=2$ in \eqv(5.lem8.2') and \eqv(5.lem8.2'bis) give \eqv(5.lem8.2)) and \eqv(5.lem8.2bis).
(Note that the bound \eqv(5.lem8.2'bis) is independent of $l$.)
The proof of Lemma \thv(5.lemma8) is complete.\endproof

\proofof{Proposition \thv(5.prop6)} By definition of an intermediate space scale, any sequence $a_n$
such that $a_n\sim b_n$ must satisfy $a_n/2^n=o(1)$. Let us first assume that $\sum_{n}a_n/2^n<\infty$.
This implies in particular that $(a_n\log n)/2^n=o(1)$ and $n/a_n=o(1)$.
Thus, using Lemma \thv(5.lemma7) with $L=2\log n$,
it follows from Borel-Cantelli Lemma that
$$
\lim_{n\rightarrow\infty}\nu_n(u,\infty)=\nu^{int}(u,\infty)\,,\text{$\P$-almost surely.}
\Eq(5.prop6.3)
$$
Together with the monotonicity of $\nu_n$ and the continuity of the limiting
function $\nu^{int}$, \eqv(5.prop6.3) entails the existence of a subset $\O^{\tau}_{1,1}\subset\O^{\tau}$
with the property that $\P(\O^{\tau}_{1,1})=1$, and such that, on $\O^{\tau}_{1,1}$,
$$
\lim_{n\rightarrow\infty}\nu_n(u,\infty)=\nu^{int}(u,\infty)\,,\quad\forall\, u>0\,.
\Eq(5.prop6.6)
$$
Similarly, using \eqv(5.lem8.2) of Lemma \thv(5.lemma8) with $L=2^n/a_n$,
it follows from
Lemma \thv(5.lemma8') and Borel-Cantelli Lemma that
$$
\eqalign{
&\lim_{n\rightarrow\infty}n\s_n^2(u,\infty)=\nu^{int}(2u,\infty)\,,\text{$\P$-almost surely.}
\cr
}
\Eq(5.prop6.4)
$$
This and the monotonicity of $\s_n^2$ allows us to conclude that
there exist a subset $\O^{\tau}_{1,2}\subset\O^{\tau}$ of full measure such that, on $\O^{\tau}_{1,2}$,
$$
\lim_{n\rightarrow\infty}n\s_n^2(u,\infty)=\nu^{int}(2u,\infty)\,,\quad\forall\, u>0\,.
\Eq(5.prop6.7)
$$
Finally, since the convergence is uniform in \eqv(5.prop6.6) then for each $0<\d\leq 1$,
on $\O^{\tau}_{1,1}$,
$
\lim_{n\rightarrow\infty}\int_0^{\d}\nu_n(u,\infty)du
=\int_0^{\d}\nu^{int}(u,\infty)du
=\d^{1-\a(\varepsilon)}\frac{\a(\varepsilon)}{1-\a(\varepsilon)}\G(\a(\varepsilon))
$.
Now $\int_0^{\d}\nu_n(u,\infty)du$ is a monotone increasing sequence having a continuous limit, and so, there exists
a subset $\O^{\tau}_{1,3}\subset\O^{\tau}$
with the property that $\P(\O^{\tau}_{1,3})=1$, such that, on $\O^{\tau}_{1,3}$,
$$
\lim_{n\rightarrow\infty}\int_0^{\d}\nu_n(u,\infty)du
=\int_0^{\d}\nu^{int}(u,\infty)du\,,\quad\forall\, 0<\d\leq 1\,.
\Eq(5.prop6.12)
$$
But this implies that, on $\O^{\tau}_{1,3}$,
$
\lim_{\d\rightarrow 0}\lim_{n\rightarrow\infty}\int_0^{\d}\nu_n(u,\infty)=0
$.
Assertion i) of the proposition now follows by taking
$\O^{\tau}_1=\O^{\tau}_{1,1}\cap\O^{\tau}_{1,2}\cap\O^{\tau}_{1,3}$.

To prove Assertion ii)  first note that by \eqv(5.lem7.1),
given $\e<1$, there exists $n_0$ such that for all $n>n_0$, for all $u>0$,
$
\nu_n(u,\infty)\leq \e+\nu^{int}(u,\infty)\leq 1 \vee 2\nu^{int}(u,\infty)
$.
Using this bound and choosing  $L=\sqrt{2^n/a_n}$ in \eqv(5.lem7.2), we obtain that for each fixed $u>0$,
$$
\P\left(\left|\nu_n(u,\infty)-\E[\nu_n(u,\infty)]\right|
\geq 2(a_n/2^n)^{1/4}\bigl(1 \vee\sqrt{2\nu^{int}(2u,\infty)}\bigr)\right)
\leq \exp\{-\sqrt{2^n/a_n}\}
\Eq(5.prop6.9)
$$
We now want to make use of Lemma 9.9 of \cite{G1}
with
$
X_n(u)=\nu_n(u,\infty)
$,
$
f_n(u)=\E[\nu_{n}(u,\infty)]
$,
$
g_n(u)\equiv g(u)=\bigl(1 \vee\sqrt{2\nu^{int}(2u,\infty)}\bigr)
$,
$
\eta_n=2(a_nL/2^n)^{1/4}$,
and
$
\rho_n=\exp\{-\sqrt{2^n/a_n}\}
$.
Indeed $\bigl(1 \vee\sqrt{2\nu^{int}(2u,\infty)}\bigr)$ is a positive decreasing function, so is  $\nu_n(u,\infty)$ for each $n$, and the properties
(9.17) of Lemma 9.9 of \cite{G1} are readily checked. Thus,
$$
\lim_{n\rightarrow 0}\P\left(
\sup_{u>0}\left\{\left|\nu_{n}(u,\infty)-\nu^{int}(u,\infty)\right|
\geq  (a_n/2^n)^{1/4}\bigl(1 \vee\sqrt{2\nu^{int}(2u,\infty)}\bigr)\right\}\right)
=0\,.
\Eq(5.prop6.10)
$$
Choosing  $L=\sqrt{2^n/a_n}$ in \eqv(5.lem8.2) of Lemma \thv(5.lemma8)
we likewise  obtain that
$$
\lim_{n\rightarrow\infty}\P\left(
\sup_{u>0}\left\{\left|n\s_n^2(u,\infty)-\nu^{int}(2u,\infty)\right|
\geq (a_n/2^n)^{1/4}\bigl(1 \vee\sqrt{2\nu^{int}(u,\infty)}\bigr)\right\}\right)
=0\,.
\Eq(5.prop6.11)
$$
To see that the last condition of Assertion ii) is satisfied
we again make use of 
 Lemma 9.9 of \cite{G1}, choosing this time
$
X_n(\d)=\int_0^{\d}\nu_n(u,\infty)du
$,
$
f_n(\d)=\int_{0}^{\d}\nu^{int}(u,\infty)du
$,
$
g_n(\d)\equiv g(\d)=\int_0^{\d}\sqrt{2\nu^{int}(2u,\infty)}du
$,
and
$
\eta_n=2(a_nL/2^n)^{1/4}
$.
Clearly,  $f_n(\d)$ and $g_n(\d)$ are positive increasing functions for each $n$,
and the leftmost relation in
(9.17)  of \cite{G1}  is satisfied, with reversed inequality,
for all $l\geq 1/\d_0$ and small enough $\d_0$.
Moreover, it follows from  \eqv(5.prop6.10)  that there exists a sequence $0<\rho_n\downarrow 0$ such that, 
setting
$$
A_n(\d)=
\left\{\left|\int_0^{\d}\nu_n(u,\infty)du-\int_{0}^{\d}\nu^{int}(u,\infty)du\right|
\geq 
2(a_n/2^n)^{1/4}\int_{0}^{\d}\sqrt{2\nu^{int}(2u,\infty)}du
\right\}\,,
\Eq(5.prop6.13)
$$
for all $n$ large enough, for all $0<\d\leq \d_0$ and small enough $\d_0$,
$
\P\left(A_n(\d)\right)\leq\rho_n
$.
Hence 
Lemma 9.9 of \cite{G1} applies, yielding
$
\lim_{n\rightarrow 0}\P\left(\sup_{0<\d<\d_0}A_n(\d)\right)=0
$.
Now by \eqv(5.lem7.1),
$
\int_0^{\d}\nu^{int}(u,\infty)du
= \d^{1-\a(\varepsilon)}\frac{\a(\varepsilon)}{1-\a(\varepsilon)}\G(\a(\varepsilon))
$
and
$
\int_0^{\d}\sqrt{\nu^{int}(u,\infty)}du
= \d^{1-\a(\varepsilon)/2}\frac{\a(\varepsilon)}{1-\a(\varepsilon)/2}\G(\a(\varepsilon))
$.
Hence we have established that there exists
$\O^{\tau}_{2,n}\subset\O^{\tau}$ with $\P(\O^{\tau}_{2,n})\geq 1-o(1)$
such that for all $n$ large enough, on $\O^{\tau}_{2,n}$, for all $0<\d\leq \d_0$ and small enough $\d_0$,
$$
\int_0^{\d}\nu_n(u,\infty)du\leq 
c_0\d^{1-\a(\varepsilon)}\sfrac{\a(\varepsilon)}{1-\a(\varepsilon)}\G(\a(\varepsilon))\,,
\Eq(5.prop6.14)
$$
for some numerical constant $0<c_0<\infty$.
This together with \eqv(5.prop6.10) and \eqv(5.prop6.11) imply Assertion ii) of the proposition.
\endproof.

\proofof{Proposition \thv(5.prop7)} By definition of a short space scale, any sequence $a_n$
such that $a_n\sim b_n$ satisfies $\sum_{n}a_n/2^n<\infty$.
Based on this observation the proof of the proposition runs along the same lines as
that of assertion (i) of Proposition \thv(5.prop6) with the following modifications.
Since on short space scales $\a(\varepsilon)=0$ it follows from
\eqv(5.lem7.1) of Lemma \thv(5.lemma7) that $\lim_{n\rightarrow\infty}\E[\nu_n(u,\infty)]=1$
(hence, the limiting function $\nu^{int}(u,\infty)$ from the statement of
Proposition \thv(5.prop6) is replaced by the constant 1). Next, $a_n$ can either satisfy
$n/a_n=o(1)$ or $a_n/n=\OO(1)$. If $n/a_n=o(1)$ the last line of \eqv(5.prop7.1)
follows from Lemma \thv(5.lemma8') and \eqv(5.lem8.2) of Lemma \thv(5.lemma8) exactly as
the last line of \eqv(5.prop6.1) follows from them.
If however $a_n/n=\OO(1)$
then
Lemma \thv(5.lemma8') yields
$\lim_{n\rightarrow\infty}a_n\E[\s_n^2(u,\infty)]=1+\lim_{n\rightarrow\infty}a_n/n$.
The last  line of \eqv(5.prop7.2)
now follows from \eqv(5.lem8.2bis) of Lemma \thv(5.lemma8), choosing e\.g\. $L=n^2$.
\endproof

\bigskip
\line{\bf 5.2 Proofs of the results of Section 1: the case of intermediate $\&$ short scales.\hfill}

In this subsection we prove the results of Section 1 that are concerned with intermediate and short scales.
These are: Assertion (i) and Assertion (ii) of Theorem \thv(1.theo1) of Subsection 1.2,
and Proposition \thv(1.prop1) and Proposition \thv(1.prop0) of Subsection 1.3.
All these results rely on the central Theorem \thv(1.theo3)  of Subsection 1.4.

We begin with results valid for intermediate scales.

\proofof{Proposition \thv(1.prop1) and Assertion (ii) of Theorem \thv(1.theo1)}
Let $r_n$ be an intermediate space scale and assume that  $\b\geq\b_c(\varepsilon)$.
Choose $\nu=\nu^{int}$ and $a_n\sim b_n$ in Conditions (A1), (A2), and (A3) (see \eqv(G1.A1) -\eqv(G1.A3)).
By the ergodic theorem of Proposition \thv(4.prop1) and the chain independent
estimates of Proposition \thv(5.prop6),  Conditions (A1), (A2),  (A3) and (A0') are
satisfied $\P$-almost surely if
$\frac{2^{m_n}}{2^n}\log n=o(1)$, and in $\P$-probability otherwise.
Thus  \eqv(G1.3.theo1.1) of Theorem \thv(1.theo3) implies that, w.r.t\. the same convergence mode as above,
$S_n(\cdot)\Rightarrow  S^{int}(\cdot)$ where $S^{int}$ is the subordinator of L\'evy measure
$\nu^{int}$. This proves Proposition \thv(1.prop1).
In addition, by \eqv(G1.3.theo1.4) of Theorem \thv(1.theo3),
$$
\lim_{n\rightarrow\infty}\CC_{n}(t,s)=\CC^{int}_{\infty}(t,s)\quad\forall\, 0\leq t<t+s\,,
\Eq(1.prop1.?)
$$
where
$
\CC^{int}_{\infty}(t,s)=\PP\left(\left\{ S^{int}(u)\,,u>0\right\}
\cap (t, t+s)=\emptyset\right)
$.
Assume first that $\b>\b_c(\varepsilon)$.
Then $S^{int}$ is a stable subordinator of index $0<\a(\varepsilon)<1$.
Thus, by the Dynkin-Lamperti Theorem in continuous time
(see e.g\. \cite{G1}, Appendix A.2, Eq\. (10.7) of Theorem 10.2),
for all $t\geq 0$ and all $\rho>0$,
$
\CC^{int}_{\infty}(t,\rho t)=\asl_{\a(\varepsilon)}(1/1+\rho)
$.
Assume next that $\b=\b_c(\varepsilon)$. Then $\a(\varepsilon)=1$, implying that
the range of $S^{int}$ is the entire positive line $[0, \infty)$.
Thus here
$
\CC^{int}_{\infty}(t,s)=\PP\left([0, \infty)\cap (t, t+s)=\emptyset\right)=0
$,
for all $0\leq t<t+s$.
Taking $s=\rho t$ in \eqv(1.prop1.?) then yields \eqv(1.theo1.1).
This proves Assertion (ii) of Theorem \thv(1.theo1).\endproof

\proofof{Proposition \thv(1.prop0) and Assertion (i) of Theorem \thv(1.theo1)}
Let $r_n$ be a short space scale and let $0<\beta<\infty$ be arbitrary.
Choose $\nu=\nu^{short}=\d_{\infty}$ and $a_n\sim b_n$ in Conditions (A1), (A2), and (A3).
We want to follow the same strategy as before and use the ergodic theorem of
Proposition \thv(4.prop1), but the latter is not useful unless $a_n\gg\t_n$.
This is why we need a lower bound on $r_n$.
Proceeding as in the proof of \eqv(2.23'') we have $\log r_n=\b\b_c(m_n/n) n(1+o(1))$.
Thus, assuming that
$
\frac 1{\b  n}\log r_n\geq 4\sqrt{{\log n}/{n}}
$
implies that $m_n\geq 4\sfrac{\log n}{\log 2}$, so that $a_n\sim b_n=2^{m_n}>n^4$.
We now easily conclude,
using Proposition \thv(4.prop1) and the estimates of Proposition \thv(5.prop7),
that Conditions  (A1), (A2), (A3) and (A0') are satisfied $\P$-almost surely.
Eq\. \eqv(G1.3.theo1.1) of Theorem \thv(1.theo3) then yields that, $\P$-almost surely,
$S_n(\cdot)\Rightarrow  S^{short}(\cdot)$ where $S^{short}$ is the
subordinator of L\'evy measure $\nu^{short}$. This proves Proposition \thv(1.prop0).
In addition, by \eqv(G1.3.theo1.4) of Theorem \thv(1.theo3),
$$
\lim_{n\rightarrow\infty}\CC_{n}(t,s)=\CC^{short}_{\infty}(t,s)\quad\forall\, 0\leq t<t+s\,,
\Eq(1.prop0.?)
$$
where
$
\CC^{short}_{\infty}(t,s)=\PP\left(\left\{ S^{short}(u)\,,u>0\right\}\cap (t, t+s)=\emptyset\right)
$.
Thus, since the range of $S^{short}$ reduces to the single point $0$,
$
\CC^{short}_{\infty}(t,s)=\PP\left(\{0\}\cap (t, t+s)=\emptyset\right)=1
$
for all $0\leq t<t+s$. This proves Assertion (i) of Theorem \thv(1.theo1)
\endproof


\bigskip


\chap{6. Extreme scales.}6


The techniques used to deal with extreme scales differ notably from those used on shorter scales.
Indeed, when $a_n\sim b_n\sim 2^n$, the convergence properties of sums such as \eqv(5.2) or \eqv(5.3)
can no longer be controlled using a classical law of large number.
The method we will use to do this, known as ``the method of common probability space'',
consists in replacing the sequence of re-scaled landscapes $(\g_n(x), x\in\VV_n)$, $n\geq 1$,
by a new sequence with identical distribution and  almost sure convergence properties.

This section closely follows Section 7 of  \cite{G1}
where this approach was first implemented.
In subsection 6.1, we give an explicit representation of the re-scaled
landscape which is valid for all extreme scales (Lemma \thv(6.lemma1))
and show that, in this representation, all random variables of interest
have an almost sure limit (Proposition \thv(6.prop4)). In subsection 6.2
we consider the model obtained by substituting the representation for
the original landscape.
For this model we state and prove the analogue of the `ergodic theorem'
of Section 4 (Proposition \thv(6.prop1)) and the associated chain independent
estimates of Section 5 (Proposition \thv(6.prop2)). Thus equipped we will be
ready, in subsection 6.3, to prove the results of Section 1 that
are concerned with extreme scales.

\bigskip
\line{\bf 6.1  A representation of the landscape.\hfill}

The representation we now introduce
is due to Lepage {\it et al\.} \cite{LWZ} and relies on an elementary property of
order statistics. We will use the following notations.
Set $N=2^n$.
Let
$\bar\t_{n}(\bar x^{(1)})\geq\dots\geq\bar\t_{n}(\bar x^{(N)})$
and
$\bar\g_{n}(\bar x^{(1)})\geq\dots\geq\bar\g_{n}(\bar x^{(N)})$
denote, respectively, the landscape and re-scaled landscape variables
$\g_n(x)=r_n^{-1}\t_n(x)$, $x\in\VV_n$, arranged in decreasing order of magnitude.
As in Section 2, set $G_n(v)=\P(\tau_n(x)>v)$, $v\geq 0$, and denote by
$
G_n^{-1}(u):=\inf\{v\geq 0 : G_n(v)\leq u\}
$, $u\geq 0$, its inverse.
Also recall that $\a=\b_c/\b$ and assume that $\b>\b_c$.

Let $(E_i, i\geq 1)$ be a sequence of i\.i\.d\. mean one exponential random variables
defined on a common probability space $(\O, \FF, \bold P)$.
For $k\geq 1$ set
$$
\eqalign{
\G_k&=\sum_{i=1}^k E_i\,,\cr
\g_k&=\G_k^{-1/\a}\,,
}
\Eq(6.4)
$$
and, for $1\leq k\leq N$, $n\geq 1$, define
$$
\g_{n}(x^{(k)})=r_n^{-1}G_n^{-1}(\G_k/\G_{N+1})\,,
\Eq(6.5)
$$
where $\{x^{(1)},\dots,x^{(N)}\}$ is a randomly chosen labelling of the $N$ elements of $\VV_n$,
all labellings being equally likely.

\lemma{\TH(6.lemma1)} {\it  For each $n\geq 1$,
$
(\bar\g_{n}(\bar x^{(1)}),\dots,\bar\g_{n}(\bar x^{(N)}))\overset d\to=(\g_{n}(x^{(1)}),\dots,\g_{n}(x^{(N)}))\,.
$
}

\proof Note that $G_n$ is non-increasing and right-continuous so that $G_n^{-1}$
is non-increasing and right-continuous. It is well known that if the random variable $U$
is a uniformly distributed on $[0,1]$ we may write $\t_n(0)\overset d\to=G_n^{-1}(U)$
(see e\.g\. \cite{Re}, page 4).
In turn it is well known (see \cite{Fe}, Section III.3) that if $(U(k), 1\leq k\leq N)$
are independent random variables uniformly distributed on $[0,1]$ then, denoting by
$\bar U_{n}(1)\leq\dots\leq \bar U_{n}(N)$ their ordered statistics,
$
(\bar U_{n}(1),\dots,\bar U_{n}(N))\overset d\to=(\G_1/\G_{N+1},\dots,\G_{N}/\G_{N+1})
$.
Combining these two facts readily yields the claim of the lemma since, by independence of
the landscape variables $\t_n(x)$, all arrangements of the $N$ variables $\G_k/\G_{N+1}$
on the $N$ vertices of $\VV_n$ are equally likely.
\endproof

Next, let $\Upsilon$  be the point process
in $M_P(\R_+)$
which has counting function
$$
\Upsilon([a,b])=\sum_{i=1}^{\infty}\1_{\{\g_k\in [a,b]\}}\,.
\Eq(6.6)
$$
\lemma{\TH(6.lemma2)} {\it $\Upsilon$ is a Poisson random measure on $(0,\infty)$
with mean measure $\mu$ given by \eqv(1.theo2.0).
}

\proof
The point process 
$
\G=\sum_{i=1}^{\infty}\1_{\{\G_k\}}
$
defines a homogeneous Poisson random measure on $[0,\infty)$ and thus,
by the mapping theorem (\cite{Re}, Proposition 3.7), setting $T(x)=x^{-1/\a}$ for $x>0$,
$\Upsilon=\sum_{i=1}^{\infty}\1_{\{T(\G_k)\}}$ is Poisson random measure on $(0,\infty)$
with mean measure $\mu(x,\infty)=T^{-1}(x)$.
\endproof

We thus established that both the ordered landscape variables and the point process $\Upsilon$
can be expressed in terms of the common sequence $(E_i, i\geq 1)$ and thus,
on the common probability space $(\O, \FF, \bold P)$. This is central to the proof of the
next proposition.

\proposition{\TH(6.prop4)} {\it Assume that $\a<1$.
Let $r_n$ be an extreme space scale. Let $f:(0,\infty)\rightarrow[0,\infty)$ be
a continuous function that obeys
$$
\int_{(0,\infty)}\min(f(u), 1)d\mu(u)<\infty\,.
\Eq(6.prop4.1)
$$
Then, $\bold P$-almost surely,
$$
\lim_{n\rightarrow\infty}\sum_{k=1}^{N}f(\g_{n}(x^{(k)}))=\sum_{k=1}^{\infty}f(\g_{k})<\infty\,.
\Eq(6.prop4.2)
$$
}

The proof of Proposition \thv(6.prop4) closely follows that of Proposition 7.3 of \cite{G1}, which itself is
strongly inspired from the proof of Proposition 3.1 of \cite{FIN}.

\proofof{Proposition \thv(6.prop4)}
Lemma \thv(2.lemma4) of Section 2 will come into use now.
By the strong law of large numbers there exists a subset $\wt\O\subset\O$
of full measure such that, for all  $n$ large enough and all  $\o\in\wt\O$,
$
\G_{N+1}=b_n(1+\l_n)
$
where
$
\l_n=o(1)
$.
From now on we assume that $\o\in\wt\O$. Thus
$$
\sum_{i=1}^{N}f(\g_{n}(x^{(i)}))
=\sum_{i=1}^{N}f\left(g_n(\G_i/(1+\l_n))\right)\,.
\Eq(6.15)
$$
where $g_n$ is defined as in \eqv(2.12).
Let us first consider the case $f(x)=x$, $x>0$.
Recall the notation $\g_i=\G^{-1/\a}_i$. For $y>0$ set
$
I(y)=\{i\geq 1 : \g_i\geq y\}
$,
$
I^c(y)=\{i\geq 1 : \g_i< y\}
$
and, for $\kappa_n=\left(b_n(1-\Phi(1/(\b\sqrt n)))\right)^{-1/\a}$,  $\d>0$,  and large enough $n$,  write:
$$
\sum_{i=1}^{N}\g_{n}(x^{(i)})=
\sum_{i\in I(\d)}\g_{n}(x^{(i)})+\sum_{i\in I(\kappa_n)\setminus I(\d)}\g_{n}(x^{(i)})+\sum_{i\in I^c(\kappa_n)}\g_{n}(x^{(i)})\,.
\Eq(6.16)
$$
From assertion (i) of Lemma \thv(2.lemma4) it follows that,
$$
\sum_{i\in I(\d)}\g_{n}(x^{(i)})\rightarrow \sum_{i\in I(\d)}\g_i\,,\quad n\rightarrow\infty\,.
\Eq(6.17)
$$
Next, by assertion (ii) of Lemma \thv(2.lemma4), for all $0<\d<1$ and some constant $0<C<\infty$, we have
$$
\sum_{i\in I(n^{-1/\a})\setminus I(\d)}\g_{n}(x^{(i)})
\leq  \sum_{i\in I(\kappa_n)\setminus I(\d)} C\G_i^{-1/\a}
=\sum_{i\in I(\kappa_n)\setminus I(\d)} C\g_i\,.
\Eq(6.18)
$$
The last sum in \eqv(6.18) is bounded above by
$
W_\d= \sum_{i:\g_i\leq \d} C\g_i
$,
and, proceeding as in (6.18)-(6.19) of \cite{G1} one gets that, choosing $d>0$ such that $\d+\a$
$
W:=\lim_{\d\rightarrow 0}W_\d=0
$
$\bold P$-almost surely.
Finally, for $i\in I^c(\kappa_n)$, $\G_i/b_n\leq 1-\Phi(1/(\b\sqrt n))$.
Being the inverse of the tail of a probability distribution, $G_n^{-1}(x)\downarrow 0$
as $x\uparrow 1$, and $G_n^{-1}(x)=0$ for all $x\geq 1$. From the calculations of the proof of
Lemma \thv(2.lemma4) we see that for small $\e>0$, $G_n^{-1}(1-\e)\approx\exp(-\b\sqrt{2n\log(1/\e)})$.
Thus, for $n$ large enough,
$$
G_n^{-1}(\G_i/[b_n(1+\l_n)])\leq G_n^{-1}([1-\Phi(1/(\b\sqrt n))]/[1+\l_n])\leq 1\,,
\Eq(6.21)
$$
so that
$$
\sum_{i\in I^c(\kappa_n)}\g_{n}(x^{(i)})
\leq \sum_{i\in I^c(\kappa_n)} r_n^{-1}
\leq b_n r_n^{-1}\,.
\Eq(6.22)
$$
Clearly $b_n r_n^{-1}\downarrow 0$ as $n\uparrow\infty$
since $b_n=2^{m_n}=2^{n(1-o(1))}$, whereas  by \eqv(2.23'), $r_n=e^{(\b\b_c n[1-o(1)])}$ where by assumption on $\b$,
$\b\b_c\geq\b_c^2=2\log 2$.
Together with \eqv(6.22) this yields,
$$
\sum_{i\in I^c(\kappa_n)}\g_{n}(x^{(i)})\rightarrow 0\,,\quad n\rightarrow\infty\,.
\Eq(6.24)
$$
Combining the previous estimates we obtain that, on a subset of $\O$ of full measure,
$$
\lim_{n\rightarrow\infty}\sum_{i=1}^{N}\g_{n}(x^{(i)})
=\lim_{\d\rightarrow 0}\sum_{i:\g_i\geq \d}\g_i
=\sum_{i=1}^{\infty}\g_{i}\,.
\Eq(6.25)
$$

The proof of \eqv(6.prop4.2) goes along the same line. We refer the reader to the
proof of (6.6) of \cite{G1} for details. This concludes the proof of Proposition \thv(6.prop4).
\endproof

\bigskip
\line{\bf 6.2 Preparations to the verification of Conditions (A1), (A2), and (A3).\hfill}

Consider the model obtained by substituting the representation
$(\g_{n}(x^{(i)}), 1\leq i\leq N)$ for the original re-scaled landscape $(\g_n(x), x\in\VV_n)$.
The aim of this subsection is to prove
the homologue, for this model, of the `ergodic theorem' (Proposition \thv(4.prop1))
and chain independent estimates (Proposition \thv(5.prop6)) of Section 4 and Section 5.

In order to distinguish the quantities $\nu_n^{J,t}(u,\infty)$, $(\s_n^{J,t})^2(u,\infty)$,
$\nu_n(u,\infty)$ and $\s_n^2(u,\infty)$, expressed in \eqv(4.7)--\eqv(4.5) in the original
landscape variables, from their expressions in the new ones , we call the latter
 $\bold v_n^{J,t}(u,\infty)$, $(\bold s_n^{J,t})^2(u,\infty)$,
$\bold v_n(u,\infty)$ and $\bold s_n^2(u,\infty)$ respectively. Their definition is otherwise unchanged.

\proposition{\TH(6.prop1)} {\it There exists a subset $\O_{0}\subset\O$ such that
$\bold P(\O_{0})=1$ and such that, on $\O_{0}$, for all large enough $n$,
the following holds: for all $t>0$ and all $u>0$,
$$
P_{\pi_n}\left(\left|\bold v_n^{J,t}(u,\infty)-E_{\pi_n}\left[\bold v_n^{J,t}(u,\infty)\right]\right|\geq\e\right)
\leq
\e^{-2}\Theta_n(t,u)\,,\quad\forall\e>0\,,
\Eq(6.prop1.1)
$$
where
$$
\eqalign{
\Theta_n(t,u)=&
\frac{k_n(t)}{a_n}\left(\frac{\bold v_n^2(u,\infty)}{2^n}+\bold s_n^2(u,\infty)\right)
+c_1\frac{\bold v_n(2u,\infty)}{ n^{2}}
\cr
+&\frac{k_n(t)}{a_n}\left(3\theta_ne^{-u/\d_n}\bold v_n(u,\infty)+\frac{2^n}{a_n}\bold v_n^2(u,\infty)e^{-c_2u}\right)\,,
}
\Eq(6.prop1.2)
$$
for some constants $0<c_1,c_2<\infty$, where $\d_n\leq n^{-\a(1+o(1))}$, and where $\theta_n$ is defined as in Proposition \thv(3.prop1).
In addition, for all $t>0$ and all $u>0$,
$$
P_{\pi_n}\left((\bold s_n^{J,t})^2(u,\infty)\geq\e'\right)\leq\frac{k_n(t)}{\e'\, a_n}\bold s_n^2(u,\infty)\,,\quad\forall\e'>0\,.
\Eq(6.prop1.3)
$$
}

\proposition{\TH(6.prop2)} {\it  Let $r_n$ be an extreme space scale and
choose $a_n\sim b_n$. Assume that  $\b\geq\b_c$ let $\nu^{ext}$ be defined in \eqv(1.prop2.2).
There exists a subset $\O_{1}\subset\O$ such that $\bold P(\O_{1})=1$ and such that, on $\O_{1}$,
the following holds: for all $u>0$,
$$
\eqalign{
\lim_{n\rightarrow\infty}\bold v_n(u,\infty)=&\nu^{ext}(u,\infty)<\infty\,,
\cr
\lim_{n\rightarrow\infty}\bold s^2_n(u,\infty)=&0\,,
}
\Eq(6.prop2.1)
$$
and
$
\lim_{\d\rightarrow 0}\lim_{n\rightarrow\infty}\int_0^{\d}\bold v_n(u,\infty)=0
$.
}

Proposition \thv(6.prop2) is a straightforward application of Proposition \thv(6.prop4) and Lemma \thv(1.lemma5)
whose proof we skip (see also \cite{G1}, (6.32)-(6.35) for a pattern of proof).

\proofof{Proposition \thv(6.prop1)} Proposition \thv(6.prop1) is a rerun of the proof
Proposition \thv(4.prop1). The only difference is in the treatment of the term \eqv(4.prop1.19).
In the new landscape variables, Lemma \eqv(4.lemma1) is not true, and its method of proof is
unadapted. To bound \eqv(4.prop1.19) we proceed as follows.
Let $T_n:=\{x^{(k)}, 1\leq k\leq n\}\subset\VV_n$ be the set of the $n$ vertices
with largest $\g_{n}(x)$. The next two lemmata collect elementary properties of $T_n$.

\lemma{\TH(6.lemma3)} {\it
There exists a subset $\O_{0,1}\subset\O$ with
$\bold P(\O_{0,1})=1$ such that, for all $\o\in\O_{0,1}$, for all large enough $n$, the following holds:
for all $x, x'\in T_n$, $x\neq x'$,
$
\dist(x,x')=\frac{n}{2}(1-\rho_n)
$
where  $\rho_n=\sqrt{\frac{8\log n}{n}}$.
}

\proof  Given $t>0$ consider the event
$
\O_{0,1}(n)=\left\{\exists_{1\leq k\neq k'\leq n}: \left|\dist\bigl(x^{(k)},x^{(k')}\bigr)-\frac{n}{2}\right|\geq t\right\}
$.
By construction, the elements of $T_n$ are drawn at random from the $\VV_n$,
independently and without replacement. Hence
$$
\textstyle
\bold P\left(\O_{0,1}(n)\right)
\leq n^2\bold P\left(\left|\dist\bigl(x^{(1)},x^{(2)}\bigr)-\frac{n}{2}\right|\geq t\right)
\sim n^2 P\left(\left|\sum_{i=1}^n\varepsilon_i-\frac{n}{2}\right|\geq t\right)\,,
\Eq(6.prop1.4)
$$
where $(\varepsilon_i, 1\leq i\leq n)$ are i.i.d\. r.v\.'s taking value 0 and 1 with probability $1/2$.
A Classical exponential Tchebychev inequality yields
$
P\left(\left|\sum_{i=1}^n\varepsilon_i-\frac{n}{2}\right|\geq t\right)\leq e^{-\frac{t^2}{2n}}
$.
Choosing $t=\sqrt{8n\log n}$, and plugging into \eqv(6.prop1.4),
$
\bold P\left(\O_{0,1}(n)\right)\leq n^{-2}
$.
Setting $\O_{0,1}=\cup_{n_0}\cap_{n>n_0}\O_{0,1}(n)$,
the claim of the lemma follows from an application of Borel-Cantelli Lemma.
\endproof

\lemma{\TH(6.lemma4)} {\it
There exists a subset $\O_{0,2}\subset\O$ with
$\bold P(\O_{0,2})=1$ such that, for all $\o\in\O_{0,2}$, for all large enough $n$,
$
\sup\{\g_{n}(x)\,,x\in\VV_n\setminus T_n\}\leq \d_n
$
where $\d_n=(1+o(1))n^{-\a(1+o(1))}$.
}

\proof Clearly
$
\sup\{\g_{n}(x)\,,x\in\VV_n\setminus T_n\}
=\sup\{\g_{n}(x^{(k)})\,,k>n\}
=\g_{n}(x^{n+1})
=r_n^{-1}G_n^{-1}\bigl(\frac{\G_{n+1}}{\G_{N+1}}\bigr)
$,
where the last equality follows from \eqv(6.5). Reasoning as in the paragraph preceding \eqv(6.15),
but applying the strong law of large numbers to both $\G_{n+1}$ and $\G_{N+1}$,
we deduce that there exists a subset $\O_{0,2}\subset\O$
of full measure such that, for all  $n$ large enough and all  $\o\in\O_{0,2}$,
$
\g_{n}(x^{n+1})
=r_n^{-1}G_n^{-1}\bigl((n/b_n)(1+\l_n)\bigr)
$.
By definition of $h_n(v)$ (see \eqv(2.13)),
$
r_n^{-1}G_n^{-1}(h_n(v))=v
$,
and by Lemma \thv(2.lemma6),
$
\g_{n}(x^{n+1})=(1+o(1))n^{-\a(1+o(1))}
$.
\endproof

We are now equipped to bound $(III)_{2,l}$. Set $\O_{0}=\O_{0,1}\cap\O_{0,2}$.
Writing $T_n^c\equiv\VV_n\setminus T_n$, and setting
$
f(y,z)=k_n(t)\pi_n(y)e^{-u[\g^{-1}_n(y)+\g^{-1}_n(z)]}p_n^{l+2}(y,z)
$,
we may decompose $(III)_{2,l}$ it into four terms,
$$
\sum_{z\in T_n^c, y\in T_n^c:y\neq z}f(y,z)
+\sum_{z\in T_n^c, y\in T_n}f(y,z)
+\sum_{z\in T_n, y\in T_n^c }f(y,z)
+\sum_{z\in T_n, y\in T_n:y\neq z}f(y,z)\,.
\Eq(6.prop1.5)
$$
To bound the first sum above we use that, by Lemma \thv(6.lemma4), for $y\in T_n^c$,
$
e^{-u[\g^{-1}_n(z)+\g^{-1}_n(y)]}\leq e^{-u/\g_n(z)}e^{-u/\d_n}
$.
Thus,
$$\eqalign{
\sum_{z\in T_n^c, y\in T_n^c:y\neq z}f(y,z)
&\leq
e^{-u/\d_n}\sum_{z\in T_n^c}
k_n(t)\pi_n(z) e^{-u/\g_n(z)}\sum_{y\in T_n^c:y\neq z}p_n^{l+2}(y,z)
\cr
&\leq
e^{-u/\d_n}\sum_{z\in T_n^c}
k_n(t)\pi_n(z) e^{-u/\g_n(z)}
\cr
&\leq
e^{-u/\d_n}\frac{k_n(t)}{a_n}\bold v_n(u,\infty)\,.
}
\Eq(6.prop1.6)
$$
The second and third sums of \eqv(6.prop1.5) are bounded just in the same way.
To deal with the last sum we use that in view of Lemma \thv(6.lemma3) the assumptions of
Proposition \thv(3.prop4) are satisfies. Consequently
$$
\eqalign{
\sum_{l=1}^{\theta_n-1}\sum_{z\in T_n, y\in T_n:y\neq z}f(y,z)
&\leq \frac{2^nk_n(t)}{a_n^2}\Bigl[a_n\sum_{z\in T_n}\pi_n(z)e^{-u/\g_n(z)}\Bigr]^2\sum_{l=1}^{\theta_n-1}p_n^{l+2}(y,z)\,,
\cr
&\leq e^{-cn}\frac{2^nk_n(t)}{a_n^2}(\bold v_n(u,\infty))^2\,,
}
\Eq(6.prop1.7)
$$
for some constant $0<c<\infty$.
Collecting \eqv(6.prop1.5), \eqv(6.prop1.6) and \eqv(6.prop1.7), and summing over $l$, we finally get,
$$
\sum_{l=1}^{\theta_n-1}(III)_{2,l}\leq 3\theta_n e^{-u/\d_n}\frac{k_n(t)}{a_n}\bold v_n(u,\infty)+
e^{-cn}\frac{2^nk_n(t)}{a_n^2}(\bold v_n(u,\infty))^2\,.
\Eq(6.prop1.8)
$$

Proposition \thv(6.prop1) is now proved just as Proposition \thv(4.prop1),
using the bound \eqv(6.prop1.8) instead of the bound  \eqv(4.lem1.2) of Lemma \eqv(4.lemma1).
\endproof

\bigskip
\line{\bf 6.3 Proofs of the results of Section 1: the case of extreme scales.\hfill}

We now prove the results of Section 1 that are concerned with extreme scales, namely:
Proposition \thv(1.prop2),
Assertion (iii) of Theorem \thv(1.theo1), Theorem \thv(1.theo2), and Lemma \thv(1.lemma5).
Again our key tool will be Theorem \thv(1.theo3) of Subsection 1.4.


We assume throughout this section that $r_n$ is an extreme space scale and that  $\b\geq\b_c$.

\proofof{Proposition \thv(1.prop2) and Assertion (iii) of Theorem \thv(1.theo1)}
Consider the model obtained by substituting the representation
$(\g_{n}(x^{(i)}), 1\leq i\leq N)$ for the original landscape $(\g_n(x), x\in\VV_n)$.
Let $\wt{\bold S}_n(\cdot)$, $\bold S_n(\cdot)$, and $\bold C_{n}(t,s)$
denote, respectively, the clock process \eqv(1.1.6), the re-scaled clock process \eqv(G1.3.2'),
and the time correlation function \eqv(1.1.8) expressed in the new landscape variables.
Choose $\nu=\nu^{ext}$ and $a_n\sim b_n$ in Conditions (A1),  (A2), and (A3)
(that is, in \eqv(G1.A1), \eqv(G1.A2), and \eqv(G1.A3), expressed of course in the new landscape variables).
By Proposition \thv(6.prop1) and Proposition \thv(6.prop2),
there exists a subset $\O_{2}\subset\O$ with
$\bold P(\O_{2})=1$, such that, on $\O_{2}$,  Conditions (A1), (A2), (A3), and (A0') are satisfied.
By  \eqv(G1.3.theo1.1) of Theorem \thv(1.theo3) we thus have that, on $\O_{2}$,
$\bold S_n(\cdot)\Rightarrow  S^{ext}(\cdot)$
where $S^{ext}$ is the (random) subordinator of L\'evy measure
$\nu^{ext}$. This proves Proposition \thv(1.prop2).

To prove Assertion (iii) of Theorem \thv(1.theo1) first note that by Lemma \thv(6.lemma1),
$$
\CC_n(t,s)\overset{d}\to=\bold C_{n}(t,s)\text{for all $n\geq 1$ and all $0\leq t<t+s$.}
\Eq(1.prop1.p4)
$$
Next, by \eqv(G1.3.theo1.4) of Theorem \thv(1.theo3) we have that, on $\O_{2}$,
$$
\lim_{n\rightarrow\infty}\bold C_{n}(t,s)=\CC^{ext}_{\infty}(t,s)\quad\forall\, 0\leq t<t+s\,,
\Eq(1.prop1.p1)
$$
where
$
\CC^{ext}_{\infty}(t,s)=\PP\left(\left\{ S^{ext}(u)\,,u>0\right\}
\cap (t, t+s)=\emptyset\right)
$.
Now, by Lemma \thv(1.lemma5), there exists a subset $\O_{3}\subset\O$ with
$\bold P(\O_{3})=1$, such that, on $\O_{3}$, $\nu^{ext}$ is
regularly varying at infinity with index $-\a$. Thus, by Dynkin-Lamperti Theorem in continuous time
applied for fixed $\o\in\O_{3}$ (see e.g\. \cite{G1}, Appendix A.2, Eq\. (10.6) of Theorem 10.2)
we get that,
$$
\lim_{t\rightarrow 0+}\CC^{ext}(t,\rho t)=\asl_{\a}(1/1+\rho)\quad\forall\,\rho>0\,.
\Eq(1.prop1.p2)
$$
Thus, by \eqv(1.prop1.p4) with  $s=\rho t$, using successively \eqv(1.prop1.p1)
and \eqv(1.prop1.p2) to pass to the limit $n\rightarrow\infty$ and $t\rightarrow 0+$,
we obtain that, for all $\rho>0$,
$
\lim_{t\rightarrow 0+}\lim_{n\rightarrow\infty}\CC^{ext}_n(t,\rho t)\overset{d}\to=\asl_{\a}(1/1+\rho)
$.
Since convergence in distribution to a constant implies convergence in probability,
the claim of Theorem \thv(1.theo1), (iii) follows.
\endproof


\proofof{Theorem \thv(1.theo2)} The proof of Theorem \thv(1.theo2)
is a re-run of the proof of Theorem 4.5 of \cite{G1} (setting $a=0$).
Note indeed that for all  $\b>\b_c$, $0\leq \a\leq 1$, which implies that
$
\int_{0}^{\infty}\nu^{ext}(u,\infty)du=\sum_{k=1}^{\infty}\g_k<\infty$
$\bold P$-almost surely.
We are thus in the realm of ``classical'' renewal theory,
in the so-called ``finite mean life time'' case. The second and first assertions of
Theorem \thv(1.theo2) will then follow, respectively, from Theorem 10.2, (ii), and
Theorem 10.4, (ii), of appendix A.2 of \cite{G1}. Their proofs use the following two
elementary observations:
$$
\frac{\int_{s}^{\infty}\nu^{ext}(u,\infty)du}{\int_{0}^{\infty}\nu^{ext}(u,\infty)du}={\CC}^{sta}_{\infty}(s)\,,\quad u>0\,,
\Eq(1.lem6.1)
$$
where ${\CC}^{sta}_{\infty}$ is defined in \eqv(1.theo2.1); Moreover, setting
$$
1-F_n(v):=\sum_{x\in\VV_n}\GG_{\a,n}(x)e^{-vc_n\l_n(x)}
=\sum_{k}\frac{\g_{n}(x^{(k)})}{\sum_{l}\g_{n}(x^{(l)})}e^{-s/\g_{n}(x^{(l)})}\,,
$$
a simple application of Proposition \thv(6.prop4) yields,
$
\lim_{n\rightarrow\infty}(1-F_n(v))=(1-F^{sta}(v)):={\CC}^{sta}_{\infty}(s)
$
$\bold P$-almost surely.
We refer the reader to \cite{G1} for details.
\endproof

It now remains to prove Lemma \thv(1.lemma5).

\proofof{Lemma \thv(1.lemma5)} To ease the notation set $\varepsilon'=1$.
Set $u^{-\a}=M$ and $f(x)=e^{-1/x}$. By \eqv(1.prop2.2) we may write
$$
u^{\a}\nu^{ext}(u,\infty)=\frac{1}{M}\sum_{k=1}^{\infty}f(M^{1/\a}\g_k)\,.
\Eq(1.lem5.3)
$$
The lemma will thus be proven if we can prove that:
$$
\lim_{M\rightarrow\infty}\frac{1}{M}\sum_{k=1}^{\infty}f(M^{1/\a}\g_k)=\a\G(\a)\text{$\bold P$-almost surely.}
\Eq(1.lem5.4)
$$
Note that it is enough for this to take the limit along
the integers since,
$f(M^{1/\a}\g_k)$ being a strictly increasing function of $M$,
$$
\frac{\lfloor M \rfloor}{M}\frac{1}{\lfloor M \rfloor}\sum_{k=1}^{\infty}f(\lfloor M \rfloor^{1/\a}\g_k)
\leq
\frac{1}{M}\sum_{k=1}^{\infty}f(M^{1/\a}\g_k)
\leq
\frac{\lceil M \rceil}{M}\frac{1}{\lceil M \rceil}\sum_{k=1}^{\infty}f(\lceil M \rceil^{1/\a}\g_k)\,.
\Eq(1.lem5.5)
$$
The claim will now follow from a classical large deviation upper bound. Set
$$
A_M=\left\{\left|
\frac{1}{M}\sum_{k=1}^{\infty}f(M^{1/\a}\g_k)-\bold E \frac{1}{M}\sum_{k=1}^{\infty}f(M^{1/\a}\g_k)
\right|\geq\d_M\right\}\,,
\Eq(1.lem5.6)
$$
where $\d_M$ is defined through
$
\d^2_M=4\a\G(\a)\frac{\log M}{M}
$.
By Tchebychev exponential inequality,
for all $\l>0$,
$$
\bold P\left(A_M\right)\leq 2 \exp\left\{-\l\d_M-\bold E ({\l}/{M})\sum_{k=1}^{\infty}f(M^{1/\a}\g_k)
+\log\bold E \exp\left\{({\l}/{M})\sum_{k=1}^{\infty}f(M^{1/\a}\g_k)\right\}
\right\}\,.
\Eq(1.lem5.7)
$$
Simple Poisson point process calculations yield
$
\bold E \frac{1}{M}\sum_{k=1}^{\infty}f(M^{1/\a}\g_k)=\a\G(\a)
$
and
$$
\log\bold E \exp\left\{({\l}/{M})\sum_{k=1}^{\infty}f(M^{1/\a}\g_k)\right\}
=-\int_{0}^{\infty}(1-e^{\frac{\l}{M} f(M^{1/\a}x)})d\mu(x)
\,.
\Eq(1.lem5.8)
$$
Furthermore, for all $k\geq 1$,
$
\int_{0}^{\infty} f^k(M^{1/\a}x)d\mu(x)=k^{-\a}M\a\G(\a)
$.
Thus
$$
-\int_{0}^{\infty}(1-e^{\frac{\l}{M} f(M^{1/\a}x)})d\mu(x)\leq \a\G(\a)\left(\l+\frac{\l^2}{4M}e^{\frac{\l}{2M}}\right)\,.
\Eq(1.lem5.9)
$$
From this last bound and the choice $\l=\d_M\frac{2M}{\a\G(\a)}$, \eqv(1.lem5.7) yields
$$
\bold P\left(A_M\right)\leq
2\exp\left\{-\frac{\d^2_MM}{\a\G(\a)}\left(2-e^{{2\d_M}/{\a\G(\a)}}\right)\right\}\leq\frac{2}{M^2}\,,
\Eq(1.lem5.10)
$$
where the last inequality follows from the definition of $\d_M$. Thus $\sum_M \bold P\left(A_M\right)\leq\infty$, and
this and the first Borel-Cantelli Lemma prove \eqv(1.lem5.4).
\endproof

\bigskip


\end